%% file: manuscript-amspreprint.tex
\pdfoutput=1
\IfFileExists{./prepreamble-amspreprint.sty}{\RequirePackage[packages,theorems,changes]{./prepreamble-amspreprint}}{}

\makeatletter
\IfFileExists{./scoop-latex/scoop-packages.sty}{\providecommand*{\input@path}{}\edef\input@path{{./scoop-latex/}\input@path}}{}
\@namedef{ver@minted.sty}{}
\@namedef{opt@minted.sty}{}
\makeatother
\PassOptionsToPackage{style = authoryear-comp, backend = biber}{biblatex}
\PassOptionsToPackage{commentmarkup = footnote}{changes}
\PassOptionsToPackage{noabbrev, capitalise, nameinlink}{cleveref}
\RequirePackage{algpseudocode}    

\documentclass[english]{amsart}

\RequirePackage{biblatex}
\RequirePackage{ltxcmds}
\IfFileExists{./preamble-amspreprint.sty}{\RequirePackage[packages,theorems,changes]{./preamble-amspreprint}}{}

\usepackage{scoop-local}

\IfFileExists{./postpreamble-amspreprint.sty}{\RequirePackage[packages,theorems,changes]{./postpreamble-amspreprint}}{}

\makeatletter
\@ifpackageloaded{changes}{
\definechangesauthor[name = {Lukas Baumgärtner}, color = {red!80!black}]{LB}
\definechangesauthor[name = {Ronny Bergmann}, color = {blue!80!black}]{RB}
\definechangesauthor[name = {Roland Herzog}, color = {orange!80!black}]{RH}
\definechangesauthor[name = {Stephan Schmidt}, color = {green!70!black}]{STS}
\definechangesauthor[name = {José Vidal-Núñez}, color = {teal}]{JVN}
\definechangesauthor[name = {Manuel Weiß}, color = {red!80!black}]{MW}
}{}
\makeatother        

\makeatletter
\@ifpackageloaded{hyperref}{%
	\hypersetup{
		pdftitle = {Mesh Denoising and Inpainting using the Total Variation of the Normal and a Shape Newton Approach},
		pdfauthor = {Lukas Baumgärtner, Ronny Bergmann, Roland Herzog, Stephan Schmidt, José Vidal-Núñez, Manuel Weiß},
		pdfkeywords = {mesh denoising, mesh inpainting, total variation of the normal vector, split Bregman iteration, shape Hessian},
		pdfcreator = {Created using the Scoop Template Engine version 1.2.0.post7+git.7a8e0c92.}
	}
}{
	\pdfinfo{
		/Title (Mesh Denoising and Inpainting using the Total Variation of the Normal and a Shape Newton Approach)
		/Author (Lukas Baumgärtner, Ronny Bergmann, Roland Herzog, Stephan Schmidt, José Vidal-Núñez, Manuel Weiß)
		/Subject ()
		/Keywords (mesh denoising, mesh inpainting, total variation of the normal vector, split Bregman iteration, shape Hessian)
		/Creator (Created using the Scoop Template Engine version 1.2.0.post7+git.7a8e0c92.)
	}
}
\makeatother

\title[Mesh Denoising and Inpainting using the TV of the Normal]{Mesh Denoising and Inpainting using the Total Variation of the Normal and a Shape Newton Approach}

\author[L. Baumgärtner]{Lukas Baumgärtner\orcidlink{0000-0003-1007-4815}}
\address[L. Baumgärtner]{Institut für Mathematik, Humboldt University of Berlin, 10099 Berlin, Germany}
\email{lukas.baumgaertner@hu-berlin.de}
\urladdr{https://www.mathematik.hu-berlin.de/en/people/mem-vz/1693318}

\author[R. Bergmann]{Ronny Bergmann\orcidlink{0000-0001-8342-7218}}
\address[R. Bergmann]{Norwegian University of Science and Technology, Department of Mathematical Sciences, NO-7041 Trondheim, Norway}
\email{ronny.bergmannn@ntnu.no}
\urladdr{https://www.ntnu.edu/employees/ronny.bergmann}

\author[R. Herzog]{Roland Herzog\orcidlink{0000-0003-2164-6575}}
\address[R. Herzog]{Interdisciplinary Center for Scientific Computing, Heidelberg University, 69120 Heidelberg, Germany}
\email{roland.herzog@iwr.uni-heidelberg.de}
\urladdr{https://scoop.iwr.uni-heidelberg.de}

\author[S. Schmidt]{Stephan Schmidt\orcidlink{0000-0002-4888-0794}}
\address[S. Schmidt]{University of Trier, Universitätsring 15, 54296 Trier, Germany}
\email{stephan.schmidt@uni-trier.de}
\urladdr{https://www.math.uni-trier.de/\string~schmidt}

\author[J. Vidal-Núñez]{José Vidal-Núñez\orcidlink{0000-0002-1190-6700}}
\address[J. Vidal-Núñez]{University of Alcalá, Department of Physics and Mathematics, 28801 Alcalá de Henares, Spain}
\email{j.vidal@uah.es}
\urladdr{https://www.uah.es/es/estudios/profesor/Jose-Vidal-Nunez/}

\author[M. Weiß]{Manuel Weiß\orcidlink{0000-0003-2164-6575}}
\address[M. Weiß]{Interdisciplinary Center for Scientific Computing, Heidelberg University, 69120 Heidelberg, Germany}
\email{roland.herzog@iwr.uni-heidelberg.de}
\urladdr{https://scoop.iwr.uni-heidelberg.de}

\thanks{This work was supported by DFG grants HE 6077/10--1 and SCHM~3248/2--1 within the Priority Program SPP~1962 (Non-smooth and Complementarity-based Distributed Parameter Systems: Simulation and Hierarchical Optimization), which is gratefully acknowledged. The work of the fifth author was partially supported by MICIIN of Spain grant AICO/2021/165.}

\date{\today}

\dedicatory{}

\begin{document}

% Insert the abstract.
\begin{abstract}
\input{abstract.tex}
\end{abstract}

% Insert the keywords.

% Insert the Mathematics Subject Classification.
\makeatletter
\ltx@ifpackageloaded{hyperref}{%
\subjclass[2010]{}
}{%
\subjclass[2010]{}
}
\makeatother

% Typeset the opening page.
\maketitle

% Insert the document body.
\input{main.tex}

% Insert the appendix.
\appendix

% Insert the bibliography.
\printbibliography

\end{document}

%% file: abstract.tex
We present a novel approach to denoising and inpainting problems for surface meshes.
The purpose of these problems is to remove noise or fill in missing parts while preserving important features such as sharp edges.
A discrete variant of the total variation of the unit normal vector field serves as a regularizing functional to achieve these goals.
In order to solve the resulting problem, we use a version of the split Bregman (ADMM) iteration adapted to the problem.
A new formulation of the total variation regularizer, as well as the use of an inexact Newton method for the shape optimization step, bring significant speed-up compared to earlier methods.
Numerical examples are included, demonstrating the performance of our algorithm with some complex 3D geometries.

%% file: main.tex
% !TEX root = manuscript-numapde-preprint.tex
\section{Introduction}
\label{section:introduction}

Meshes are widely employed in computer graphics and computer vision, where they are utilized to represent generated shapes and real geometries.
Meshes are produced by 3D~scanners and can be efficiently processed numerically with appropriate software.
On the other hand, the scanning process for geometry acquisition leads to unavoidable errors in the form of noise, or even missing parts.
The task of removing such noise while preserving relevant features is known as \emph{mesh denoising}.
When missing parts of the geometry must be filled in, we speak of \emph{mesh inpainting}, also known as \emph{hole filling}.
The aforementioned problems have been of interest to the \emph{image processing} community since the late 1980s; see for instance \cite{CasellesChambolleNovaga:2015:1} for an overview.
Mesh denoising and related problems have a similarly rich history, and many algorithms exist.
We refer the reader, \eg, to \cite{BotschPaulyKobbeltAlliezLevyBischoffRoessl:2007:1,ChenWeiWang:2022:1} for surveys.

In this paper, we present an approach for simultaneous denoising and inpainting of triangular surface meshes, using the total variation (TV) of the cell-wise constant normal vector as regularizer.
While our definition of the TV of the normal agrees with the definition previously used in \cite{WuZhengCaiFu:2015:1} and \cite{BergmannHerrmannHerzogSchmidtVidalNunez:2020:2}, our algorithmic treatment of an ensuing variational problem is different.
The authors of \cite{WuZhengCaiFu:2015:1} approximated the TV of the normal for algorithmic purposes and used an augmented Lagrangian approach to solve the minimization problem.
The advance over \cite{BergmannHerrmannHerzogSchmidtVidalNunez:2020:2} in this paper is two-fold.
First, we use a reformulation of the total variation of the normal functional, which allows us to simplify the numerical algorithm.
More precisely, the parallel transport of the Lagrange multipliers between tangent spaces of the $2$-sphere, a relatively expensive step, can be dropped.
Second, instead of using a first-order method for the optimization with respect to the nodal coordinates of the mesh, we present a significantly faster, second-order globalized Newton method, based on second-order shape calculus.

The rest of the paper is organized as follows.
\Cref{section:Related_work} reviews related work in mesh denoising.
\Cref{section:Directional_TV_of_the_normal} contains background on the notion of total variation of the normal vector field.
We also present our novel, equivalent formulation of this functional there, that is favorable for numerical purposes.
\Cref{section:numerical_realization} is devoted to the derivation of the split Bregman iteration (ADMM) for general problems involving the total variation of the normal as regularizer.
The main emphasis is on the derivation and presentation of Newton's method for the shape optimization subproblem in \cref{subsection:Newton_method}.
In \cref{section:Mesh_Denoising_Problem} and \cref{section:Mesh_Inpainting_Problem} we discuss the application of the split Bregman algorithm specifically for mesh denoising and inpainting problems, respectively, and present numerical results.

\subsection*{Notation}

Throughout the paper, we denote vector-valued quantities such as the normal vector~$\bn$ by bold-face symbols.
The Euclidean inner product between such vectors $\bn_1, \bn_2$ in $\R^3$ is denoted by $\bn_1 \cdot \bn_2$ or $\bn_1^\transp \bn_2$.
The symbol $\id$ denotes the identity matrix.

\section{Related Work on Mesh Denoising}
\label{section:Related_work}

We begin our review of related work by recalling that \emph{surface fairing} (or \emph{surface smoothing}) should be clearly distinguished from mesh denoising.
The goal of the former is to smooth a given geometry.
By contrast, the emphasis of the latter is to remove spurious information from the geometry while preserving sharp features such as edges, characterized by a significant change in the normal vector.

Various types of mesh denoising approaches exist in the literature.
Methods based on diffusion can be classified as either isotropic or anisotropic.
They can often be traced back to perimeter minimization and curvature flows.
On the one hand, isotropic diffusion methods are applied for mesh fairing and their main characteristic is that they do not take into account the geometric features.
Examples of such methods are \cite{Field:1988:1,Taubin:1995:1,DesbrunMeyerSchroederBarr:1999:1} and \cite{VollmerMenclMueller:1999:1}.
These methods generally have the drawback that they may lead to surface shrinkage and tend to blur geometric features.
On the other hand, anisotropic diffusion methods can give good recovery results, but meshes with sharp edges can still be a challenge for them.
To overcome this drawback, several approaches have been developed in the literature, as detailed in the following.

\subsection*{Anisotropic Diffusion}
\label{subsection:Anisotropic_Diffusion}

Anisotropic diffusion methods take into account feature directions while filtering the normal vector.
The references \cite{BajajXu:2003:1,ClarenzDiewaldRumpf:2000:1,TasdizenWhitakerBurchardOsher:2002:1} and \cite{HildebrandtPolthier:2004:1} fall into this category.
Generally speaking, these methods can preserve genuine features such as edges, but they may be numerically unstable during the diffusion process \cite{ZhaoQinZengXuDong:2018:1}.

\subsection*{Bilateral Filtering}
\label{subsection:Bilateral_Filtering}

Bilateral filtering approaches are one-stage iterative methods, that use the normal and tangential directions to the surface to determine a filter at every vertex of the mesh.
They then compute a displacement correction for the vertex and update its position.
This estimation, however, may be inaccurate due to the presence of noise in the mesh, which may lead to edges not being well preserved.
In this category, we mention \cite{FleishmanDroriCohenOr:2003:1,JonesDurandDesbrun:2003:1}, where the authors extend the bilateral filter method in imaging processing from \cite{TomasiManduchi:1998:1} to denoise 3D meshes.

\subsection*{Normal Filtering and Vertex Update}
\label{subsection:Normal_Filtering_and_Vertex_Update}

The class of normal filtering and vertex update methods are characterized as two-stage methods, that first filter the facet normals and subsequently update the vertex positions according to the filtered normals.
Nevertheless, small-scale features might get blurred since most approaches filter the facet normals by averaging over neighboring normals.
Methods in this class differ \wrt the treatment of the facet normal vector, while the vertex update usually is straightforward.
We mention \cite{OhtakeBelyaevSeidel:2002:1,SunRosinMartinLangbein:2007:1,ZhengFuAuTai:2011:1,ZhangDengZhangBouazizLiu:2015:1,WangFuLiuTongLiuGuo:2015:1} and more recently, \cite{WangLiuTong:2016:1,YadavReitebuchPolthier:2018:1} and \cite{CentinSignoroni:2018:1}.

\subsection*{\texorpdfstring{$L_0$-Minimization}{L0-Minimization}}
\label{subsection:L0_Minimization}

Another popular class of mesh denoising methods is based on the so-called $L_0$-minimization, which combines vertex and normal regularization; see \cite{HeSchaefer:2013:1} and \cite{ZhaoQinZengXuDong:2018:1}.
The $L_0$-term introduces sparsity into a discrete gradient operator describing the variation of the surface.
One of the weak points of this class of methods is that the non-convexity of the model can lead to a high demand in computational resources.

\subsection*{Vertex Classification and Denoising}
\label{subsection:Vertex_Classification_and_Denoising}

Vertex classification and denoising methods can also be understood as two-stage methods, that first classify the vertices of the mesh and subsequently apply a denoising method per class or cluster of vertices.
As examples of such methods we mention \cite{WeiYuPangWangQinLiuHeng:2015:1,LuDengChen:2016:1,WangZhangYu:2012:1,ZhuWeiYuWangQinHeng:2013:1,WeiLiangPangWangLiWu:2017:1}.
Although the idea of separating vertices in homogeneous classes seems to be promising, the presence of noise can make the vertex classification difficult or unreliable; see \cite{LuDengChen:2016:1}.
As a consequence, this class of methods depends heavily on the level of noise present in the mesh.

\subsection*{Mumford-Shah}
\label{subsection:Mumford-Shah}

The famous Mumford-Shah functional, introduced in \cite{MumfordShah:1989:1}, is very similar in its effect to the total variation regularizer we utilize.
Both are able to reconstruct discontinuities in the optimization variable.
The Mumford-Shah functional achieves this by identifying areas without discontinuities, and it applies smoothing there while ignoring the discontinuities between such areas.
Originally designed for segmentation and denoising of images, the Mumford-Shah functional, or rather the Ambrosio-Tortorell approximation, see \cite{AmbrosioTortorelli:1990:1}, was recently also used for mesh denoinsing and inpainting in \cite{BonneelCoeurjollyGuethLachaud:2018:1}.

\subsection*{Total Variation for Mesh Denoising}
\label{subsection:Total_Variation_for_Mesh_Denoising}

As with many of the approaches discussed above, the use of the total variation (TV) seminorm for denoising originates from imaging, see~\cite{RudinOsherFatemi:1992:1,ChanEsedogluParkYip:2006:1,ChambolleCasellesCremersNovagaPock:2010:1}.
It is well-known that the TV-seminorm can preserve high frequency features in images.
Despite this advantage, solutions to total variation based image denoising problems generally suffer from the so-called staircasing effect.
Moreover, the non-differentiability of the TV-seminorm requires customized algorithms; see for instance \cite{CasellesChambolleNovaga:2015:1} and the references therein.
To the best of our knowledge, \cite{TasdizenWhitakerBurchardOsher:2002:1} was the first to mention the total variation functional for mesh optimization.
\cite{ElseyEsedoglu:2009:1} then proposed an analogue of the ROF model based on the total variation of the Gaussian curvature.

In this paper, we use the concept of total variation of the normal vector that was proposed in \cite{WuZhengCaiFu:2015:1}.
It arises from the definition of the total variation for manifold valued data, see \cite{LellmannStrekalovskiyKoetterCremers:2013:1}, because the unit normal vector belongs to the $2$-sphere.
We previously used the TV of the normal as regularizer for mesh denoising problems in \cite{BergmannHerrmannHerzogSchmidtVidalNunez:2019:3,HerrmannHerzogSchmidtVidalNunez:2022:1}.
Moreover, we employed it for a geometric inverse problem subject to an electrical impedance tomography problem in \cite{BergmannHerrmannHerzogSchmidtVidalNunez:2020:2}.

In contrast to our previous works, we propose in the present paper a novel, equivalent formulation of the total variation of the normal functional that is favorable for numerical purposes.
In contrast to \cite{WuZhengCaiFu:2015:1}, we do not need to approximate the TV of the normal.

\section{Total Variation of the Normal Vector}
\label{section:Directional_TV_of_the_normal}

For the remainder of the paper, suppose that $\Gamma$ is a triangulated surface embedded in $\R^3$ composed of flat triangles~$T$, edges~$E$ and vertices~$\vertex$.
The collections of all triangles, edges and vertices are denoted by $\cT$, $\cE$ and $\cV$, respectively.
We recall that \cite{WuZhengCaiFu:2015:1} proposed the \emph{total absolute edge-lengthed supplementary angle of the dihedral angle} (TESA)
\begin{equation}\label{eq:Discrete_TV_angle}
	\sum_{E \in \cE} \theta_E \, \abs{E}_2
\end{equation}
as a regularizing functional for the purpose of surface denoising.
Here $\abs{E}_2$ represents the length of the edge $E$ and $\theta_E$ is the exterior dihedral angle between two neighboring triangles, \ie, the angle between the respective outward normal vectors.
In \cite{BergmannHerrmannHerzogSchmidtVidalNunez:2020:2}, we introduced the \emph{total variation of the normal}
\begin{equation}\label{eq:Discrete_TV_log}
	\sum_{E \in \cE} \abs{\logarithm{\bn_+}{\bn_-}}_2 \, \abs{E}_2
\end{equation}
and proved that it coincides with \eqref{eq:Discrete_TV_angle}.
Here the \enquote{$+$} and \enquote{$-$} refer to the two sides of an edge~$E$ and $\bn_+$, $\bn_-$ are the outer unit normal vectors of the respective triangles on either side of~$E$.
The normal vectors belong to the unit $2$-sphere~$\Sphere^2$ in $\R^3$.
On this manifold, the logarithmic map $\logarithm{\bn_+}{\bn_-}$ returns the tangent vector pointing from $\bn_+$ to $\bn_-$.
This map is well-defined whenever $\bn_+ \neq -\bn_-$ and it reads
\begin{equation}
	\label{eq:logarithm}
	\logarithm{\bn_+}{\bn_-}
	=
	\begin{cases}
		0
		,
		&
		\text{ if } \bn_+ = \bn_-
		\\
		\displaystyle \arccos \paren(){\bn_+ \cdot \bn_-} \, \frac{\bn_- - (\bn_+ \cdot \bn_-) \, \bn_+}{\abs{\bn_- - (\bn_+ \cdot \bn_-) \, \bn_+}}
		,
		&
		\text{ otherwise}
		.
	\end{cases}
\end{equation}
When we view the tangent space~$\tangent{\bn_+}[\Sphere^2]$ at $\bn_+$ as a subspace of $\R^3$, the Euclidean length $\abs{\logarithm{\bn_+}{\bn_-}}_2$ equals the geodesic distance in $\Sphere^2$ (the angle) between $\bn_+$ and $\bn_-$.
Note that the functional \eqref{eq:Discrete_TV_log} is non-differentiable whenever at least two neighboring normals $\bn_\pm$ agree.

\begin{remark}
	\label{remark:log_mu_parallel}
	For $\bn_+ \neq \bn_-$, the vector
	\begin{equation}
		\label{eq:normalized_log}
		\frac{\bn_- - (\bn_+ \cdot \bn_-) \, \bn_+}{\abs{\bn_- - (\bn_+ \cdot \bn_-) \, \bn_+}}
	\end{equation}
	and thus the logarithmic map $\logarithm{\bn_+}{\bn_-}$ is parallel to the so-called co-normal vector $\bmu_+$.
	The latter is the unit vector uniquely defined by the following properties: it lies in the same plane as the triangle $T_+$, it is perpendicular to the shared edge~$E$ and points outward from $T_+$; see \cref{figure:Co-Normal}.
	Furthermore, the relation
	\begin{equation}
	\label{remark:log_mu_parallel2}
		\frac{\bn_- - (\bn_+ \cdot \bn_-) \, \bn_+}{\abs{\bn_- - (\bn_+ \cdot \bn_-) \, \bn_+}}
		=
		\sign \paren(){\bmu_+ \cdot \bn_-} \, \bmu_+
	\end{equation}
	holds.
\end{remark}

\tdplotsetmaincoords{75}{140}
\begin{figure}[htb]
	\centering
	\begin{subfigure}{0.49\linewidth}
		\centering
		\begin{tikzpicture}[tdplot_main_coords]
			% Variables
			\pgfmathsetmacro{\rotationangleC}{0}
			\pgfmathsetmacro{\rotationangleD}{130}
			\pgfmathsetmacro{\arrowscale}{0.3}

			% Starting point from which to rotate
			\pgfmathsetmacro{\Sx}{1.5}
			\pgfmathsetmacro{\Sy}{4}
			\pgfmathsetmacro{\Sz}{0}

			% Define vertices of the shared edge
			\pgfmathsetmacro{\Ax}{0}
			\pgfmathsetmacro{\Ay}{0}
			\pgfmathsetmacro{\Az}{0}
			\pgfmathsetmacro{\Bx}{3}
			\pgfmathsetmacro{\By}{0}
			\pgfmathsetmacro{\Bz}{0}

			% Rotate the starting points forward and backward
			\pgfmathsetmacro{\Cx}{\Sx}
			\pgfmathsetmacro{\Cy}{{cos(\rotationangleC)*\Sy - sin(\rotationangleC)*\Sz}}
			\pgfmathsetmacro{\Cz}{{sin(\rotationangleC)*\Sy + cos(\rotationangleC)*\Sz}}
			\pgfmathsetmacro{\Dx}{\Sx}
			\pgfmathsetmacro{\Dy}{{cos(\rotationangleD)*\Sy - sin(\rotationangleD)*\Sz}}
			\pgfmathsetmacro{\Dz}{{sin(\rotationangleD)*\Sy + cos(\rotationangleD)*\Sz}}

			\pgfmathsetmacro{\ABx}{{\Bx-\Ax}}
			\pgfmathsetmacro{\ABy}{{\By-\Ay}}
			\pgfmathsetmacro{\ABz}{{\Bz-\Az}}
			\pgfmathsetmacro{\ACx}{{\Cx-\Ax}}
			\pgfmathsetmacro{\ACy}{{\Cy-\Ay}}
			\pgfmathsetmacro{\ACz}{{\Cz-\Az}}
			\pgfmathsetmacro{\ADx}{{\Dx-\Ax}}
			\pgfmathsetmacro{\ADy}{{\Dy-\Ay}}
			\pgfmathsetmacro{\ADz}{{\Dz-\Az}}

			% Define shared edge vertices
			\coordinate (A) at (\Ax,\Ay,\Az);
			\coordinate (B) at (\Bx,\By,\Bz);
			\coordinate (M) at  ($(A)!0.5!(B)$);

			% Define standalone vertices
			\coordinate (C) at (\Cx,\Cy,\Cz);
			\coordinate (D) at (\Dx,\Dy,\Dz);

			% Compute the normal wrt the C-triangle
			% Define normal n_+
			\tdplotcrossprod(\ABx,\ABy,\ABz)(\ACx,\ACy,\ACz)
			\pgfmathsetmacro{\nplusx}{\arrowscale*\tdplotresx}
			\pgfmathsetmacro{\nplusy}{\arrowscale*\tdplotresy}
			\pgfmathsetmacro{\nplusz}{\arrowscale*\tdplotresz}

			% Compute n_-
			\tdplotcrossprod(\ADx,\ADy,\ADz)(\ABx,\ABy,\ABz)
			\pgfmathsetmacro{\nminusx}{\arrowscale*\tdplotresx}
			\pgfmathsetmacro{\nminusy}{\arrowscale*\tdplotresy}
			\pgfmathsetmacro{\nminusz}{\arrowscale*\tdplotresz}

			% Compute co-normal mu_+
			\pgfmathsetmacro{\muplusx}{\nplusx}
			\pgfmathsetmacro{\muplusy}{{cos(90)*\nplusy - sin(90)*\nplusz}}
			\pgfmathsetmacro{\muplusz}{{sin(90)*\nplusy + cos(90)*\nplusz}}

			% Compute co-normal mu_-
			\pgfmathsetmacro{\muminusx}{\nminusx}
			\pgfmathsetmacro{\muminusy}{{cos(-90)*\nminusy - sin(-90)*\nplusz}}
			\pgfmathsetmacro{\muminusz}{{sin(-90)*\nminusy + cos(-90)*\nplusz}}

			% Draw co-normals
			\draw [TolVibrantBlue, ->, very thick] (M) -- ++ (\muplusx,\muplusy,\muplusz) node [anchor = south]{$\bmu_+$};
			\draw [TolVibrantBlue, ->, very thick] (M) -- ++ (\muminusx,\muminusy,\muminusz) node [anchor = west]{$\bmu_-$};

			% Draw the triangles
			\draw [thick,fill = TolVibrantBlue, opacity = 0.5] (A) -- (B) -- (C) -- cycle;
			\draw [thick,fill = TolVibrantBlue, opacity = 0.5] (A) -- (B) -- (D) -- cycle;

			% Draw the normals
			\draw [TolVibrantRed, ->, very thick] (M) -- ++ (\nplusx,\nplusy,\nplusz) node [anchor = south](N){$\bn_+$};
			\draw [TolVibrantRed, ->, very thick] (M) -- ++ (\nminusx,\nminusy,\nminusz) node [anchor = west]{$\bn_-$};

			% Draw log
			\pgfmathsetmacro{\innernmu}{\nplusx*\muminusx+\nplusy*\muminusy+\nplusz*\muminusz }
			\pgfmathsetmacro{\signnmu}{sign(\innernmu)}

			\pgfmathsetmacro{\innernn}{\nplusx*\nminusx+\nplusy*\nminusy+\nplusz*\nminusz }
			\pgfmathsetmacro{\normnplus}{sqrt(\nplusx*\nplusx+\nplusy*\nplusy+\nplusz*\nplusz) }
			\pgfmathsetmacro{\normnminus}{sqrt(\nminusx*\nminusx+\nminusy*\nminusy+\nminusz*\nminusz) }

			\pgfmathsetmacro{\angle}{acos(\innernn/\normnplus/\normnminus)}
			\pgfmathsetmacro{\arccos}{acos(\innernn/\normnplus/\normnminus)/360*2*3.1415}

			\coordinate (N) at ($(M)+(\nplusx,\nplusy,\nplusz)$);
			\draw [black,->, very thick] (N) -- ++ (\arccos*\signnmu*\muplusx,\arccos*\signnmu*\muplusy,\arccos*\signnmu*\muplusz) node [anchor = west]{$\logarithm{\bn_+}{\bn_-}$};

			\tdplotsetrotatedcoordsorigin{(M)}
			\path [fill = black, fill opacity = 0.2, tdplot_rotated_coords] (M) -- plot [variable = \t, domain = {\rotationangleC+90}:{\rotationangleD-90}](xyz spherical cs:radius = \normnplus, longitude = 0.001, latitude = \t);
		\end{tikzpicture}
	\end{subfigure}
	\hfill
	\begin{subfigure}{0.49\linewidth}
		\begin{tikzpicture}[tdplot_main_coords]
			% Starting point from which to rotate
			\pgfmathsetmacro{\Sx}{1.5}
			\pgfmathsetmacro{\Sy}{4}
			\pgfmathsetmacro{\Sz}{0}

			% Variables
			\pgfmathsetmacro{\rotationangleC}{-50}
			\pgfmathsetmacro{\rotationangleD}{180}
			\pgfmathsetmacro{\arrowscale}{0.3}

			% Starting point from which to rotate
			% Define vertices of the shared edge
			\pgfmathsetmacro{\Ax}{0}
			\pgfmathsetmacro{\Ay}{0}
			\pgfmathsetmacro{\Az}{0}
			\pgfmathsetmacro{\Bx}{3}
			\pgfmathsetmacro{\By}{0}
			\pgfmathsetmacro{\Bz}{0}

			% Rotate the starting points forward and backward
			\pgfmathsetmacro{\Cx}{\Sx}
			\pgfmathsetmacro{\Cy}{{cos(\rotationangleC)*\Sy - sin(\rotationangleC)*\Sz}}
			\pgfmathsetmacro{\Cz}{{sin(\rotationangleC)*\Sy + cos(\rotationangleC)*\Sz}}
			\pgfmathsetmacro{\Dx}{\Sx}
			\pgfmathsetmacro{\Dy}{{cos(\rotationangleD)*\Sy - sin(\rotationangleD)*\Sz}}
			\pgfmathsetmacro{\Dz}{{sin(\rotationangleD)*\Sy + cos(\rotationangleD)*\Sz}}

			\pgfmathsetmacro{\ABx}{{\Bx-\Ax}}
			\pgfmathsetmacro{\ABy}{{\By-\Ay}}
			\pgfmathsetmacro{\ABz}{{\Bz-\Az}}
			\pgfmathsetmacro{\ACx}{{\Cx-\Ax}}
			\pgfmathsetmacro{\ACy}{{\Cy-\Ay}}
			\pgfmathsetmacro{\ACz}{{\Cz-\Az}}
			\pgfmathsetmacro{\ADx}{{\Dx-\Ax}}
			\pgfmathsetmacro{\ADy}{{\Dy-\Ay}}
			\pgfmathsetmacro{\ADz}{{\Dz-\Az}}

			% Define shared edge vertices
			\coordinate (A) at (\Ax,\Ay,\Az);
			\coordinate (B) at (\Bx,\By,\Bz);
			\coordinate (M) at  ($(A)!0.5!(B)$);

			% Define standalone vertices
			\coordinate (C) at (\Cx,\Cy,\Cz);

			\coordinate (D) at (\Dx,\Dy,\Dz);

			% Compute the normal wrt the C-triangle
			% Define normal n_+
			\tdplotcrossprod(\ABx,\ABy,\ABz)(\ACx,\ACy,\ACz)
			\pgfmathsetmacro{\nplusx}{\arrowscale*\tdplotresx}
			\pgfmathsetmacro{\nplusy}{\arrowscale*\tdplotresy}
			\pgfmathsetmacro{\nplusz}{\arrowscale*\tdplotresz}

			% Compute n_-
			\tdplotcrossprod(\ADx,\ADy,\ADz)(\ABx,\ABy,\ABz)
			\pgfmathsetmacro{\nminusx}{\arrowscale*\tdplotresx}
			\pgfmathsetmacro{\nminusy}{\arrowscale*\tdplotresy}
			\pgfmathsetmacro{\nminusz}{\arrowscale*\tdplotresz}

			% Compute co-normal mu_+
			\pgfmathsetmacro{\muplusx}{\nplusx}
			\pgfmathsetmacro{\muplusy}{{cos(90)*\nplusy - sin(90)*\nplusz}}
			\pgfmathsetmacro{\muplusz}{{sin(90)*\nplusy + cos(90)*\nplusz}}

			% Compute co-normal mu_-
			\pgfmathsetmacro{\muminusx}{\nminusx}
			\pgfmathsetmacro{\muminusy}{{cos(-90)*\nminusy - sin(-90)*\nplusz}}
			\pgfmathsetmacro{\muminusz}{{sin(-90)*\nminusy + cos(-90)*\nplusz}}

			% Draw the triangles
			\draw [thick,fill = TolVibrantBlue, opacity = 0.5] (A) -- (B) -- (C) -- cycle;
			\draw [thick,fill = TolVibrantBlue, opacity = 0.5] (A) -- (B) -- (D) -- cycle;

			% Draw co-normals
			\draw [TolVibrantBlue, ->, very thick] (M) -- ++ (\muplusx,\muplusy,\muplusz) node [anchor = south]{$\bmu_+$};
			\draw [TolVibrantBlue, ->, very thick] (M) -- ++ (\muminusx,\muminusy,\muminusz) node [anchor = south]{$\bmu_-$};

			% Draw the normals
			\draw [TolVibrantRed, ->, very thick] (M) -- ++ (\nplusx,\nplusy,\nplusz) node [anchor = west](N){$\bn_+$};
			\draw [TolVibrantRed, ->, very thick] (M) -- ++ (\nminusx,\nminusy,\nminusz) node [anchor = south]{$\bn_-$};

			% Draw log
			\pgfmathsetmacro{\innernmu}{\nplusx*\muminusx+\nplusy*\muminusy+\nplusz*\muminusz }
			\pgfmathsetmacro{\signnmu}{sign(\innernmu)}

			\pgfmathsetmacro{\innernn}{\nplusx*\nminusx+\nplusy*\nminusy+\nplusz*\nminusz }
			\pgfmathsetmacro{\normnplus}{sqrt(\nplusx*\nplusx+\nplusy*\nplusy+\nplusz*\nplusz) }
			\pgfmathsetmacro{\normnminus}{sqrt(\nminusx*\nminusx+\nminusy*\nminusy+\nminusz*\nminusz) }

			\pgfmathsetmacro{\angle}{acos(\innernn/\normnplus/\normnminus)}
			\pgfmathsetmacro{\arccos}{acos(\innernn/\normnplus/\normnminus)/360*2*3.1415}

			\coordinate (N) at ($(M)+(\nplusx,\nplusy,\nplusz)$);
			\draw [black, ->, very thick] (N) -- ++ (\arccos*\signnmu*\muplusx,\arccos*\signnmu*\muplusy,\arccos*\signnmu*\muplusz) node [anchor = west]{$\logarithm{\bn_+}{\bn_-}$};

			\tdplotsetrotatedcoordsorigin{(M)}
			\path [fill = black, fill opacity = 0.2, tdplot_rotated_coords] (M) -- plot [variable = \t, domain = {\rotationangleC+90}:{\rotationangleD-90}](xyz spherical cs:radius = \normnplus, longitude = 0.001, latitude = \t);
		\end{tikzpicture}
	\end{subfigure}
	\caption{Illustration of the geodesic distance (angle) between normals $\bn_+$ and $\bn_-$ and the logarithmic map $\logarithm{\bn_+}{\bn_-}$ of two triangles $T_+$, $T_-$, which share the edge $E$. The triangles' co-normals are $\bmu_+$ and $\bmu_-$. Note that $\logarithm{\bn_+}{\bn_-}$ is parallel to $\bmu_+$.}
	\label{figure:Co-Normal}
\end{figure}

In this paper, we propose another reformulation of \eqref{eq:Discrete_TV_log}, which gives rise to simpler algorithms.
It is based on the observation
\begin{equation}
	\label{eq:abs_logarithm}
	\abs{\logarithm{\bn_+}{\bn_-}}_2
	=
	\abs{\arccos \paren(){\bn_+ \cdot \bn_-}}
	\quad
	\text{by \eqref{eq:logarithm}}
	.
\end{equation}
	On the other hand, we have
	\begin{alignat}{2}
		(\logarithm{\bn_+}{\bn_-}) \cdot \bmu_+
		&
		=
		\arccos \paren(){\bn_+ \cdot \bn_-} \, \frac{\bn_- - (\bn_+ \cdot \bn_-) \, \bn_+}{\abs{\bn_- - (\bn_+ \cdot \bn_-) \, \bn_+}} \cdot \bmu_+
		&
		&
		\quad
		\text{by \eqref{eq:logarithm}}
		\notag
		\\
		&
		=
		\arccos \paren(){\bn_+ \cdot \bn_-} \, \sign \paren(){\bmu_+ \cdot \bn_-}
		&
		&
		\quad
		\text{by \eqref{remark:log_mu_parallel2}}
		.
		\label{eq:reformulation_of_geodesic_distance2}
	\end{alignat}
Together, \eqref{eq:abs_logarithm} and \eqref{eq:reformulation_of_geodesic_distance2} show
\begin{equation}
	\label{eq:reformulation_of_geodesic_distance}
	\abs{\logarithm{\bn_+}{\bn_-}}_2
	=
	\abs{(\logarithm{\bn_+}{\bn_-}) \cdot \bmu_+}
	=
	\arccos \paren(){\bn_+ \cdot \bn_-}
	.
\end{equation}
We thus obtain the following re-formulation of \eqref{eq:Discrete_TV_log}:
\begin{equation}\label{eq:Discrete_TV_log_novel}
	\sum_{E \in \cE} \abs{E}_2 \, \abs{(\logarithm{\bn_+}{\bn_-}) \cdot \bmu_+}
	.
\end{equation}

Compared to the common definition of the discrete total variation semi-norm in imaging, which involves the absolute value of the difference of function values on neighboring cells, the term $\abs{(\logarithm{\bn_+}{\bn_-}) \cdot \bmu_+}$ in \eqref{eq:Discrete_TV_log_novel} appears to be highly nonlinear.
However, it agrees with the geodesic distance in \eqref{eq:Discrete_TV_log} and is thus the natural extension of the absolute value of the difference for $\Sphere^2$-valued data.
In addition, \eqref{eq:reformulation_of_geodesic_distance2} can be viewed as the signed distance between neighboring normal vectors, and $\sign \paren(){\bmu_+ \cdot \bn_-}$ distinguishes the convex from the concave situation; see \cref{figure:Co-Normal}.

\section{Numerical Realization via a Split Bregman Algorithm}
\label{section:numerical_realization}

In the following, we consider problems of the type
\begin{equation}
	\label{eq:general_problem_formulation}
	\text{Minimize}
	\quad
	\cF(\Gamma)
	+
	\beta \sum_{E \in \cE} \abs{E}_2 \, \abs{(\logarithm{\bn_+}{\bn_-}) \cdot \bmu_+}
	\quad
	\text{\wrt\ }
	\Gamma
	,
\end{equation}
which feature the new formulation \eqref{eq:Discrete_TV_log_novel} of the total variation of the normal as a regularizer.
The optimization variable in \eqref{eq:general_problem_formulation} is the surface mesh $\Gamma$, that has fixed connectivity but unknown vertex positions.
The feasible region can be shown to be an open submanifold of $\R^{3 \, \numberofvertices}$, where $\numberofvertices$ is the number of vertices, using techniques similar to \cite[Proposition~3.4]{HerzogLoayzaRomero:2022:1}.
Note that all of $\abs{E}_2$, $\bn_+$, $\bn_-$ as well as $\bmu_+$ depend on $\Gamma$.
We omit this dependency for brevity in our notation.
The regularization parameter $\beta$ is assumed to be positive.

In this section, we allow the shape functional $\cF$ to be generic but smooth \wrt variations of the vertex positions of $\Gamma$.
Typical choices for $\cF$ are fidelity or tracking terms and, depending on its support, \eqref{eq:general_problem_formulation} can be a regularized shape optimization, denoising or an inpainting problem.

As was mentioned in the introduction, we consider an alternating direction method of multipliers (ADMM) to deal with the non-smooth terms $\abs{(\logarithm{\bn_+}{\bn_-}) \cdot \bmu_+}$.
Due to the similarity of the regularization term with an $\ell_1$-norm, the ADMM takes the form of a split Bregman algorithm~\cite{GoldsteinOsher:2009:1}.
To derive the method, we introduce an additional variable $d$ on the skeleton of the mesh.
On each edge~$E$, the value $d_E$ is scalar and constant and we impose the constraint $d_E = (\logarithm{\bn_+}{\bn_-}) \cdot \bmu_+$.
This leads to the following re-formulation of \eqref{eq:general_problem_formulation}:
\begin{equation}
	\label{eq:general_problem_formulation_with_constraints}
	\begin{aligned}
		\text{Minimize}
		\quad
		&
		\cF(\Gamma)
		+
		\beta \sum_{E \in \cE} \abs{E}_2 \, \abs{d_E}
		\quad
		\text{\wrt\ }
		(\Gamma, d)
		\\
		\text{subject to}
		\quad
		&
		d_E
		=
		(\logarithm{\bn_+}{\bn_-}) \cdot \bmu_+
		\text{ on each edge }
		E
		.
	\end{aligned}
\end{equation}
The Lagrange multiplier vector $\lambda$ associated with the constraints is identified with a constant scalar function $\lambda_E$ on each edge $E$.
Given the augmentation parameter $\rho > 0$, the augmented Lagrangian for \eqref{eq:general_problem_formulation_with_constraints} reads
\begin{multline*}
	\widetilde{\cL}_\rho(\Gamma, d, \lambda)
	\coloneqq
	\cF(\Gamma)
	+
	\beta \sum_{E \in \cE} \abs{E}_2 \, \abs{d_E}
	\\
	+
	\sum_{E \in \cE} \abs{E}_2 \lambda_E \paren[big][]{(\logarithm{\bn_+}{\bn_-}) \cdot \bmu_+ - d_E}
	+
	\frac{\rho}{2} \sum_{E \in \cE} \abs{E}_2 \paren[big][]{(\logarithm{\bn_+}{\bn_-}) \cdot \bmu_+ - d_E}^2
	.
\end{multline*}
Substituting $b = \tfrac{\lambda}{\rho}$ and completing the squares, we obtain the equivalent form
\begin{multline}
	\label{eq:augmented_Lagrangian}
	\cL_\rho(\Gamma, d, b)
	\coloneqq
	\cF(\Gamma)
	+
	\beta \sum_{E \in \cE} \abs{E}_2 \, \abs{d_E}
	\\
	+
	\frac{\rho}{2} \sum_{E \in \cE} \abs{E}_2 \paren[big][]{d_E - (\logarithm{\bn_+}{\bn_-}) \cdot \bmu_+ - b_E}^2
	-
	\frac{\rho}{2} \sum_{E \in \cE} \abs{E}_2 b_E^2
	.
\end{multline}
Split Bregman iterations are characterized by the fact that the augmented Lagrangian $\cL_\rho$ is minimized independently with respect to the vertex positions of $\Gamma$ and auxiliary variable~$d$, respectively.
Notice that in other settings, the last term in \eqref{eq:augmented_Lagrangian} is usually omitted since it is independent of the primary variables.
Here, however, it depends on the vertex positions of~$\Gamma$ through the edge lengths, and therefore we keep the term.

At the end of each iteration, an update of the rescaled Lagrange multiplier~$b$ is performed according to $b_E \gets b_E + (\logarithm{\bn_+}{\bn_-}) \cdot \bmu_+ - d_E$.

Even though $\cL_\rho$ is non-smooth in the $d$-variable, the minimization \wrt $d$ is manageable since the problem completely decouples edge-by-edge.
Moreover, the problem for each edge has a well-known closed-form solution in terms of the soft thresholding or shrinkage operator:
\begin{equation}\label{eq:Classical_shrinkage_for_d_E_non_directional}
	d_E
	=
	\max \paren[auto]\{\}{\abs{w_E} - \beta \rho^{-1}, \; 0} \sign{(w_E)}
	,
\end{equation}
where $w_E \coloneqq (\logarithm{\bn_+}{\bn_-}) \cdot \bmu_+ + b_E$.

To minimize $\cL_\rho$ \wrt the vertex positions of the surface mesh $\Gamma$, we use a variational approach based on shape calculus, that is discussed in~\cref{subsection:shape_derivatives}.
Independently of the choice of the functional $\cF$, said minimization is challenging since $\abs{E}_2$, $\logarithm{\bn_+}{\bn_-}$ and $\bmu_+$ depend on the vertex positions in a nonlinear way.
The dependency is, however, smooth, as long as the vertex positions remain in the feasible open set, \ie, as long as there are no degenerate triangles nor self-intersection.
Therefore, we can employ a globalized Newton scheme in the same variational manner.
Since the Newton method only solves one smooth subproblem of a split Bregman iteration, a limited number of Newton updates are sufficient.

We summarize the proposed split Bregman method for the generic problem \eqref{eq:general_problem_formulation} in \cref{algorithm:split_Bregman}.
Further details, in particular concerning the globalized Newton scheme and the adaptive selection of the penalty parameter~$\rho$, will be given in \cref{subsection:Newton_method,subsection:penalty_parameter_selection} respectively.

\begin{algorithm}[htb]
	\caption{Split Bregman method for problem class \eqref{eq:general_problem_formulation}.}
	\label{algorithm:split_Bregman}
	\begin{algorithmic}[1]
		\Require initial guess $\sequence{\Gamma}{0}$ for the surface mesh
		\Require regularization parameter $\beta > 0$
		\Require initial penalty parameter $\rho > 0$
		\Ensure approximate solution of \eqref{eq:general_problem_formulation}
		\State Set $\sequence{b}{0} \coloneqq 0$ and $\sequence{d}{0} \coloneqq 0$
		\State Set $k \coloneqq 0$
		\While{not converged}
		\State Set $\sequence{d}{k+1} \coloneqq \argmin \cL_\rho(\sequence{\Gamma}{k}, d, \sequence{b}{k})$ with respect to $d$; see \eqref{eq:Classical_shrinkage_for_d_E_non_directional}
		\label{line:Split_Bregman_discrete_tv_step}
		\State Perform at most $\iter_\textup{max}$~Newton steps to obtain $\sequence{\Gamma}{k+1}$, an approximate minimizer of $\cL_\rho(\Gamma, \sequence{d}{k+1}, \sequence{b}{k})$ with respect to the vertex positions in $\Gamma$, starting from the initial guess~$\sequence{\Gamma}{k}$; see \cref{subsection:Newton_method} and \cref{algorithm:Newton_solve}
		\label{line:Split_Bregman_discrete_shape_step}
		\State Set $\sequence{b_E}{k+1} \coloneqq \sequence{b_E}{k} + (\logarithm{\bn_+^{k+1}}{\bn_-^{k+1}}) \cdot \bmu_+^{k+1} - \sequence{d_E}{k+1}$ for all edges $E$
		\label{line:Split_Bregman_discrete_multiplier_step}
		\State Adapt the penalty parameter $\rho$; see \cref{subsection:penalty_parameter_selection}
		\label{line:Split_Bregman_penalty_parameter_selection}
		\State Set $k \coloneqq k+1$
		\EndWhile
	\end{algorithmic}
\end{algorithm}

Let us discuss the main differences to the Riemannian split Bregman method proposed in \cite[Algorithm~1]{BergmannHerrmannHerzogSchmidtVidalNunez:2020:2}.
\begin{enumeratearabic}
	\item
		We employ here a globalized Newton scheme for the inexact minimization of $\cL_\rho$ with respect to the vertex positions, instead of a gradient method.
		This leads to much faster convergence of \cref{algorithm:split_Bregman} but requires the evaluation of discrete second-order shape derivatives.

	\item
		We use here the novel formulation \eqref{eq:Discrete_TV_log_novel} of the total variation of the normal regularizer instead of the original formulation \eqref{eq:Discrete_TV_log}.
		This allows us to choose the auxiliary variables $d_E = (\logarithm{\bn_+}{\bn_-}) \cdot \bmu_+$ as scalars, with associated scalar Lagrange multipliers $b_E$.
		In \cite{BergmannHerrmannHerzogSchmidtVidalNunez:2020:2}, by contrast,  we used tangent-space valued auxiliary variables
		\begin{equation}
			\label{eq:contraint_riemannian_admm}
			\widehat \bd_E
			=
			\logarithm{\bn_+}{\bn_-} \in \tangent{\bn_+}[\Sphere^2]
		\end{equation}
		and Lagrange multipliers $\widehat \bb_E$ which also belonged to $\tangent{\bn_+}[\Sphere^2]$.

		Although the minimization of the augmented Lagrangian with respect to $\widehat \bd_E$ can be expressed explicitly in terms of a (vectorial) shrinkage operation in the respective tangent space, the multiplier estimates needed to be parallely transported once per iteration into the new tangent space $\tangent{\bn_+}[\Sphere^2]$, following the update of the geometry and the implied update of the normal vectors.
		We refer the reader to \cite[Algorithm~1]{BergmannHerrmannHerzogSchmidtVidalNunez:2020:2} for details.
		By contrast, our new \cref{algorithm:split_Bregman} is significantly simpler and less expensive in that regard because it involves only scalar shrinkage operations for $d_E$ and, more importantly, does not require the parallel transport step for the Lagrange multiplier estimates.

	\item
		Another, albeit minor change compared to \cite{BergmannHerrmannHerzogSchmidtVidalNunez:2020:2} is the order of minimization steps \wrt the geometry and $d$ in each iteration of the algorithm.
		Here we first minimize \wrt $d$ and subsequently \wrt the vertex positions of $\Gamma$, instead of the other way around.
		This order turned out to yield faster convergence in practice.
\end{enumeratearabic}

\subsection{First- and Second-Order Shape Derivatives of the Augmented Lagrangian}
\label{subsection:shape_derivatives}

\Cref{line:Split_Bregman_discrete_shape_step} of the proposed split Bregman scheme (\cref{algorithm:split_Bregman}) requires the (approximate) solution of a smooth optimization problem.
As one of the novelties of the present work, we use Newton's method for this discrete shape optimization problem.
The evaluation of the first- and second-order derivatives of all quantities $\abs{E}_2$, $\bn_+$, $\bn_-$ and $\bmu_+$ that depend on the vertex positions is nontrivial.
It is tempting to use algorithmic differentiation (AD) for this purpose.
Unfortunately, as seen in \eqref{eq:logarithm} and \eqref{eq:reformulation_of_geodesic_distance}, geometric insight is needed to prevent the occurrence of singular fractions.
General AD tools are unaware of these transformations and produce numerically unstable code for the problem at hand.
In particular, computing the total derivative of $\cL_\rho$ with respect to all vertex positions of $\Gamma$ by such means is numerically unstable in the frequent situation that there are co-planar triangles, \ie, $\bn_+ = \bn_-$.
Therefore, we use shape calculus to obtain expressions that can be evaluated in a stable way.

In order to define shape derivatives in our discrete setting, we consider deformation fields $\bV$ on the surface~$\Gamma$ that are continuous ($C$) and piecewise linear ($\cP_1$) on the triangles of the mesh, \ie,
\begin{equation*}
	\bV \in \CG{1}(\Gamma, \R^3)
	\coloneqq
	\setDef[auto]{\bV \in C(\Gamma, \R^3)}{\restr{\bV}{T} \in \cP_1(T, \R^3)}
	.
\end{equation*}
For sufficiently small $\varepsilon > 0$, any such field defines a family of perturbed surface meshes
\begin{equation}
	\label{eq:pertubation_of_identity}
	\Gamma_\varepsilon[\bV]
	\coloneqq
	\setDef[auto]{x + \varepsilon \, \bV(x)}{x \in \Gamma}
\end{equation}
with the same connectivity.
The first-order shape derivative of a scalar function $J(\Gamma)$ in the direction of $\bV$ is defined by
\begin{equation}
	\label{eq:shape_derivative}
	\d J(\Gamma)[\bV]
	\coloneqq
	\lim_{\varepsilon \searrow 0} \frac{J(\Gamma_\varepsilon) - J(\Gamma)}{\varepsilon}
\end{equation}
and hence $\d J(\Gamma)[\bV] = 0$ for all $\bV \in \CG{1}(\Gamma, \R^3)$ is a first-order necessary optimality condition.
Notice that in the literature, shape derivatives in the form of \eqref{eq:shape_derivative} are also called material derivatives.
We refer the reader to \cite{DelfourZolesio:1992:1} and \cite[Section~2.23, Section~3.9]{SokolowskiZolesio:1992:1} as well as \cite{HamMitchellPaganiniWechsung:2018:1} for a thorough background.

Notice that \eqref{eq:shape_derivative} readily generalizes to vector-valued quantities such as $\d \bn(\Gamma)[\bV]$.
Furthermore, the regular chain rule applies.
Therefore, to compute $\d \cL(\Gamma, d, b)[V]$ for fixed $d, b$, we investigate the shape derivatives of the quantities $\bn$ and $\abs{E}_2$.
To this end, we require the tangential derivative $D_\Gamma \bV$ on $\Gamma$, that is well-defined in the interior of each triangle.
While $D_\Gamma \bV$ can be defined intrinsically, we view it, triangle-by-triangle, as $D_\Gamma \bV = D\widetilde{\bV} - D\widetilde{\bV} \bn \bn^\transp$, where $\widetilde{\bV}$ is any smooth extension of $\bV$ to a neighborhood of the triangle, and $D\widetilde{\bV}$ is its classical Jacobian.
For the normal vector~$\bn$, we have the identity
\begin{equation}
	\label{eq:simple_shape_derivatives:1}
	\d \bn(\Gamma)[\bV]
	=
	- (D_\Gamma \bV)^\transp \bn
	\quad
	\text{in the interior of each triangle}
	.
\end{equation}
This result can be found, for instance, in~\cite[Eq.~(3.168)]{SokolowskiZolesio:1992:1} or \cite[Lemma~3.3.6]{Schmidt:2010:1}.

By extending the classical shape derivative for manifolds of co-dimension $1$, \ie, surfaces in $\R^3$, from~\cite[Eq.~(2.172)]{SokolowskiZolesio:1992:1}, to manifolds of co-dimension $2$, \ie, lines in $\R^3$, we obtain
\begin{equation}
	\label{eq:simple_shape_derivatives:3}
	\d \abs{E}_2(\Gamma)[\bV]
	=
	\d \paren[bigg](){\int\limits_E 1 \ \d\ell}[\bV]
	=
	\int\limits_E \div_E\bV \d \ell
	=
	\abs{E}_2 \div_E \bV
	\quad
	\text{on $E$}
	,
\end{equation}
where $\div_E \bV$ is the constant edge-intrinsic divergence.
Utilizing $\bt_E$, a unit vector tangent to $E$ with arbitrary but fixed orientation, it can be defined via $\div_E \bV = \bt_E^\transp (D_\Gamma \bV) \, \bt_E$.

The Jacobian $D_\Gamma \bV$ of $\bV \in \CG{1}(\Gamma, \R^3)$ is constant per triangle and discontinuous.
Nevertheless, $(D_\Gamma \bV) \, \bt_E$ is well-defined on $E$ and independent of which of the two adjacent values of $D_\Gamma \bV$ is used.

We can now proceed to evaluate the first-order shape derivative of the augmented Lagrangian \eqref{eq:augmented_Lagrangian}.
Notice that $d_E$ and $b_E$ have no dependency on the mesh, they are just constant values that get transported along with the mesh.
Hence, their material derivative is zero.
Using the product and chain rules and in particular \eqref{eq:simple_shape_derivatives:3}, we obtain
\begin{multline}
	\label{derivative_augmented_Lagrangian}
	\d \cL_\rho(\Gamma, d, b)[\bV]
	=
	\d \cF(\Gamma)[\bV]
	+
	\beta \sum_{E \in \cE} \abs{E}_2 \, \abs{d_E} \, \div_E \bV
	\\
	+
	\frac{\rho}{2} \sum_{E \in \cE} \abs{E}_2 \paren[big][]{d_E - (\logarithm{\bn_+}{\bn_-}) \cdot \bmu_+ - b_E}^2 \div_E \bV
	-
	\frac{\rho}{2} \sum_{E \in \cE} \abs{E}_2 \, b_E^2 \div_E \bV
	\\
	+
	\rho \sum_{E \in \cE} \abs{E}_2 \paren[big][]{d_E - (\logarithm{\bn_+}{\bn_-}) \cdot \bmu_+ - b_E} \d \paren[big](){-(\logarithm{\bn_+}{\bn_-}) \cdot \bmu_+}[\bV]
	.
\end{multline}

It remains to compute the shape derivative of the term $(\logarithm{\bn_+}{\bn_-}) \cdot \bmu_+$.
To this end, we rewrite the term according to \eqref{eq:reformulation_of_geodesic_distance2}, \ie, we have
\begin{equation}
	\label{eq:reformulation_of_geodesic_distance3}
	(\logarithm{\bn_+}{\bn_-}) \cdot \bmu_+
	=
	\sign \paren(){\bmu_+ \cdot \bn_-} \arccos \paren[auto](){\bn_+ \cdot \bn_-}
	.
\end{equation}
We proceed by the chain rule and evaluate the partial derivatives of \eqref{eq:reformulation_of_geodesic_distance3} \wrt $\bn_+$ and  $\bn_-$.
We exclude the co-planar case $\bn_+ = \bn_-$ for now and observe that $\bn_+ = -\bn_-$ would result in a non-conforming mesh.
The exclusion of this case ensures that $\bmu_+ \cdot \bn_-$ has constant sign for small enough perturbations $\varepsilon\bV$.
Consequently, the derivative of $\sign \paren(){\bmu_+ \cdot \bn_-}$ \wrt $\bn_-$ and $\bmu_+$ is zero and its dependency on the vertex positions can be ignored when differentiating \eqref{eq:reformulation_of_geodesic_distance3}.

We now seek to compute the tangent derivative $D_{\Sphere^2}\arccos \paren[auto](){\bn \cdot \bn_-}$ along $\Sphere^2$.
To this end, we emphasize that $\bn$ is the spatial coordinate of $\Sphere^2$.
Furthermore, when embedding $\Sphere^2$ into $\R^3$, the outer normal is also again $\bn$.
Hence, we have
\begin{equation}
	\label{eq:TangentDerivativeNormalS2}
	\begin{aligned}
		D_{\Sphere^2}\arccos \paren[auto](){\bn \cdot \bn_-} &= D\arccos \paren[auto](){\bn \cdot \bn_-} -D\paren[auto](){\arccos \paren[auto](){\bn \cdot \bn_-}}\bn\bn^\transp
		\\
		&
		=
		D\arccos \paren[auto](){\bn \cdot \bn_-}(\id - \bn\bn^\transp)
		\\
		&
		=
		\frac{-\bn_-}{\sqrt{1-(\bn \cdot \bn_-)^2}}(\id - \bn\bn^\transp) = \frac{-(\bn_- - (\bn \cdot \bn_-) \, \bn)}{\sqrt{1-(\bn \cdot \bn_-)^2}}
		.
	\end{aligned}
\end{equation}
A simple calculation shows that
\begin{align*}
	\sqrt{1-(\bn_+ \cdot \bn_-)^2} = \abs{\bn_- - (\bn_+ \cdot \bn_-) \, \bn_+}
\end{align*}
holds.
Hence, the resulting vector is normalized and parallel to $\bmu_+$; see \cref{remark:log_mu_parallel}.
Altogether we obtain
\begin{subequations}
	\label{eq:derivative2}
	\begin{align}
		\MoveEqLeft
		\restr[Big]{D_{\Sphere^2} \paren[auto][]{\sign \paren(){\bmu_+ \cdot \bn_-} \arccos \paren[auto](){\bn \cdot \bn_-}}}{\bn = \bn_+}
		\notag
		\\
		&
		=
    -
		\sign \paren(){\bmu_+ \cdot \bn_-} \frac{\bn_- - (\bn_+ \cdot \bn_-) \, \bn_+}{\abs{\bn_- - (\bn_+ \cdot \bn_-) \, \bn_+}}
		=
		- \bmu_+
		.
		\label{eq:derivative2:1}
		\intertext{Proceeding in the same way for the derivative with respect to $\bn_-$, we obtain}
		\MoveEqLeft
		\restr[Big]{D_{\Sphere^2} \paren[auto][]{\sign \paren(){\bmu_+ \cdot \bn}\arccos \paren[auto](){\bn \cdot \bn_+}}}{\bn = \bn_-}
		\notag
		\\
		&
		=
		-\sign \paren(){\bmu_+ \cdot \bn_-} \frac{\bn_+ - (\bn_+ \cdot \bn_-) \, \bn_-}{\abs{\bn_+ - (\bn_+ \cdot \bn_-) \, \bn_-}}
		=
		- \bmu_-
		.
		\label{eq:derivative2:2}
	\end{align}
\end{subequations}
From here we also infer that the assumption $\bn_+ \neq \bn_-$ is no longer necessary, since both expressions in \eqref{eq:derivative2} are continuous even at $\bn_+ = \bn_-$.
Equation~\eqref{eq:reformulation_of_geodesic_distance3}, the chain rule on $\Sphere^2$ and \eqref{eq:derivative2} yield
\begin{align}
	\MoveEqLeft
	- \d \paren[big](){ (\logarithm{\bn_+}{\bn_-}) \cdot \bmu_+ }[\bV]
	\notag
	\\
	&
	=
	-
	\d \paren[big][]{\sign \paren(){\bmu_+ \cdot \bn_-} \arccos \paren[auto](){\bn_+ \cdot \bn_-}}[\bV]
	\quad
	\text{by \eqref{eq:reformulation_of_geodesic_distance3}}
	\notag
	\\
	&
	=
	- \restr[Big]{D_{\Sphere^2} \paren[auto][]{\sign \paren(){\bmu_+ \cdot \bn_-} \arccos \paren[auto](){\bn \cdot \bn_-}}}{\bn = \bn_+}
	\d\bn_+[\bV]
	\notag
	\\
	&
	\quad
	- \restr[Big]{D_{\Sphere^2} \paren[auto][]{\sign \paren(){\bmu_+ \cdot \bn}\arccos \paren[auto](){\bn \cdot \bn_+}}}{\bn = \bn_-}
	\d \bn_-[\bV]
	\notag
	\\
	&
	=
	\bmu_+
	\cdot
	\d
	\bn_+[\bV]
	+
	\bmu_-
	\cdot
	\d
	\bn_-
	[\bV]
	\quad
	\text{by \eqref{eq:derivative2}}
	.
	\label{eq:shape_derivative_log}
\end{align}
This clarifies the last term in \eqref{derivative_augmented_Lagrangian}, since the shape derivatives $\d \bn_+[\bV]$ and $\d \bn_-[\bV]$ are given by extending \eqref{eq:simple_shape_derivatives:1} from the respective triangle to the edge.
Notice that the intermediate terms in \eqref{eq:shape_derivative_log} cannot be stably evaluated by tools of algorithmic differentiation since \eqref{eq:TangentDerivativeNormalS2} yields $0/0$ whenever $\bn_+ = \bn_-$ holds.
This is due to the fact that the $\arccos$ function is not differentiable at~$1$.

The closed form expression \eqref{eq:shape_derivative_log} also allows for a rather simple computation of the second-order shape derivative $\d^2 \cL_\rho(\Gamma, d, b)[\bV,\bW]$ in the directions $\bV, \bW \in \CG{1}(\Gamma, \R^3)$.

Reapplying the same differentiation rules to a shape derivative expression in the sense of
\begin{equation}
	\label{eq:second-order-shape-derivative:naive}
	\d \paren[auto](){\d J[\bV]} [\bW]
	=
	\restr[auto]{\frac{\d}{\d \varepsilon_2}}{\varepsilon_2 = 0} \restr[auto]{\frac{\d}{\d \varepsilon_1}}{\varepsilon_1 = 0} J\paren[big](){\setDef{x + \varepsilon_1\bW(x) + \varepsilon_2\bV(x + \varepsilon_1\bW(x))}{x \in \Gamma}}
\end{equation}
is problematic, a behavior not known from directional derivatives in $\R^n$.
Expression \eqref{eq:second-order-shape-derivative:naive} generally leads to a non-symmetric form different from the Hessian.
A correct analytic expression that coincides with the discrete derivative needed for computations is obtained by computing
\begin{equation}
	\label{eq:second-order-shape-derivative:corrected}
  \d^2 J[\bV, \bW]
  \coloneqq
	\restr[auto]{\frac{\d}{\d \varepsilon_2}}{\varepsilon_2 = 0} \restr[auto]{\frac{\d}{\d \varepsilon_1}}{\varepsilon_1 = 0} J\paren[big](){\setDef{x + \varepsilon_1 \bV(x) + \varepsilon_2 \bW(x)}{x \in \Gamma}}
  .
\end{equation}
The symmetric expression $\d^2 J[\bV, \bW]$ can be obtained from $\d \paren[auto](){\d J[\bV]} [\bW]$ when $\d \bV[\bW] = 0$ is assumed.
For more details, we refer the reader, \eg, to \cite{SchmidtSchulz:2023:1}.
Note that in our discrete setting, the property $\d \bV[\bW] = 0$ is automatically fulfilled for arbitrary functions $\bV, \bW \in \CG{1}(\Gamma, \R^3)$.
This is because $\d^2 J[\bV, \bW]$ can be interpreted as a deformation in regard to some arbitrary fixed reference domain, where $\bV$ is not transported by $\bW$.
Here, the finite element reference element has the role of a fixed reference geometry, see also \cite{HamMitchellPaganiniWechsung:2018:1}.
Alternatively, the initial guess of the shape optimization problem can be used as reference; see \cite{SchmidtSchulz:2023:1}.
The use of starshaped domains in~\cite{EpplerHarbrechtSchneider:2007:1} is the forerunner of this concept.

To evaluate the shape Hessian of the Lagrangian $\d^2 \cL_\rho[\bV, \bW]$, we require the second-order shape derivative of the normal vector $\d^2\bn$, the derivative of the co-normal $\d\bmu$ on an edge, and the second-order shape derivative of $\d^2\abs{E}_2$.
Starting with the normal $\bn$, we obtain from \cite[Eq.~(22)]{SchmidtGraesserSchmid:2022:1}
\begin{multline}
	\label{eq:dn2}
	\d^2 \bn(\Gamma)
	[\bV,\bW]
	=
	\paren[auto](){D_\Gamma \bV}^\transp \paren[auto](){D_\Gamma \bW}^\transp \bn
	+
	\paren[auto](){D_\Gamma \bW}^\transp \paren[auto](){D_\Gamma \bV}^\transp \bn
	\\
	-
	\paren[big](){\bn^\transp \paren[auto](){D_\Gamma \bW} \paren[auto](){D_\Gamma \bV}^\transp \bn} \bn
\end{multline}
on the interior of each triangle.

We proceed with the shape derivative of the co-normal $\bmu$ on the edge~$E$ of a particular triangle~$T$.
Using again $\bt_E$, a unit vector tangent to $E$ with arbitrary but fixed orientation, we have
\begin{equation}
	\label{eq:simple_shape_derivatives:2}
	\d \bmu(\Gamma)[\bV]
	=
	(D_\Gamma \bV) \, \bmu
	-
	\paren(){\bmu^\transp (D_\Gamma \bV) \bmu} \, \bmu
	-
	\paren[big](){\bt_E^\transp \paren[big][]{D_\Gamma \bV + (D_\Gamma \bV)^\transp} \bmu} \, \bt_E
\end{equation}
on~$E$.
This result follows from \cite[Theorem~3.2]{BergmannHerrmannHerzogSchmidtVidalNunez:2020:1}, when we choose the orthonormal basis of the triangle's tangent space as $\bxi_1 = \bt_E$ and $\bxi_2 = \bmu$.

Furthermore, the second material derivative of the edge length is needed, that is, to the best of our knowledge, not in the literature yet.
Notice that the following computation generalizes to arbitrary curve integrals.
\begin{lemma}
	\label{lemma:second-derivative-of-edge-length}
	The second material derivative of the edge length $\abs{E}_2$ is given by
	\begin{equation}
		\label{eq:d2E}
		\d^2 \abs{E}_2[\bV,\bW]
		=
		\abs{E}_2 \, \bt_E^\transp \, (D_\Gamma\bW)^\transp (\id - \bt_E \bt_E^\transp ) (D_\Gamma\bV) \, \bt_E
		.
	\end{equation}
\end{lemma}
\begin{proof}
	As in the case of the Lagrangian, we differentiate the edge length twice in the sense of \eqref{eq:second-order-shape-derivative:naive} but assume $\d\bV[\bW] = 0$.
	In this process we will use that
	\begin{equation*}
		\d \bt_E [\bV]
		=
		\paren[auto](){D_\Gamma\bV} \bt_E - \bt_E \paren[auto](){\bt_E^\transp \paren[auto](){D_\Gamma\bV} \bt_E)}
	\end{equation*}
	from \cite[Theorem~3.2]{BergmannHerrmannHerzogSchmidtVidalNunez:2020:1}, as well as
	\begin{equation*}
		\d(D_\Gamma\bV)[\bW]
		=
		- \paren[auto](){D_\Gamma\bV} \paren[auto](){D_\Gamma\bW} + D_\Gamma\paren[auto](){\d\bV[\bW]}
		=
		- \paren[auto](){D_\Gamma\bV} \paren[auto](){D_\Gamma\bW}
	\end{equation*}
	from \cite[eq.(25)]{Berggren:2010:1}.

	We are now in the position to evaluate $\d^2 \abs{E}_2[\bV,\bW] = \d(\d \abs{E}_2[\bV])[\bW]$ via
	\begin{align*}
		\MoveEqLeft
		\d(\d \abs{E}_2[\bV])[\bW]
		\\
		&
		=
		\d(\abs{E}_2 \paren[auto](){\bt_E^\transp \paren[auto](){D_\Gamma\bV} \bt_E}[\bW]
		\\
		&
		=
		\d(\abs{E}_2)[\bW]\paren[auto](){\bt_E^\transp \paren[auto](){D_\Gamma \bV} \bt_E}
		+
		\abs{E}_2 \d \paren[auto](){\bt_E^\transp \paren[auto](){D_\Gamma\bV} \bt_E}[\bW]
		\\
		&
		=
		\abs{E}_2 \paren[auto](){\bt_E^\transp \paren[auto](){D_\Gamma\bW} \bt_E} \paren[auto](){\bt_E^\transp \paren[auto](){D_\Gamma\bV} \, \bt_E}
		+
		\abs{E}_2 \paren[big](){\bt_E^\transp \d\paren[auto](){D_\Gamma\bV}[\bW] \, \bt_E}
		\\
		&
		\quad
		+
		\abs{E}_2 \paren[auto](){\d\bt_E[W]^\transp \paren[auto](){ \paren[auto](){D_\Gamma\bV} + \paren[auto](){D_\Gamma\bV}^\transp} \bt_E}
		\\
		&
		=
		\abs{E}_2 \paren[auto](){\bt_E^\transp \paren[auto](){D_\Gamma\bW} \bt_E} \paren[auto](){\bt_E^\transp \paren[auto](){D_\Gamma\bV} \bt_E}
		-
		\abs{E}_2 \paren[big](){\bt_E^\transp \paren[big](){D_\Gamma\bV} \paren[auto](){D_\Gamma\bW}\bt_E}
		\\
		&
		\quad
		+
		\abs{E}_2 \paren[big](){\bt_E^\transp \paren[big](){D_\Gamma\bW}^\transp \paren[auto](){ \paren[auto](){D_\Gamma \bV} + \paren[auto](){D_\Gamma\bV}^\transp } \bt_E}
		\\
		&
		\quad
		-
		\abs{E}_2 \paren[auto](){\bt_E^\transp \paren[auto](){D_\Gamma\bW}^\transp \bt_E} \paren[auto](){\bt_E^\transp \paren[auto](){\paren[auto](){D_\Gamma\bV} + \paren[auto](){D_\Gamma\bV}^\transp } \bt_E}
		\\
		&
		=
		\abs{E}_2 \paren[big](){\bt_E^\transp \paren[auto](){D_\Gamma\bW}^\transp \paren[auto](){D_\Gamma\bV} \bt_E}
		-
		\abs{E}_2 \paren[big](){\bt_E^\transp \paren[auto](){D_\Gamma\bW}^\transp \bt_E \bt_E^\transp \paren[auto](){D_\Gamma\bV} \bt_E}
		\\
		&
		=
		\abs{E}_2 \paren[auto](){\bt_E^\transp \paren[auto](){D_\Gamma\bW}^\transp \paren[auto](){\id - \bt_E \bt_E^\transp } \paren[auto](){D_\Gamma \bV} \bt_E}
		.
	\end{align*}
\end{proof}

Finally, the second-order shape derivative of $\logarithm{\bn_+}{\bn_-}$ obtained by repeated shape differentiation \eqref{eq:second-order-shape-derivative:naive} reads
\begin{align}
	\MoveEqLeft
	\d^2 \paren[big](){-(\logarithm{\bn_+}{\bn_-}) \cdot \bmu_+}[\bV, \bW]
	\notag
	\\
	&
	=
	\d \paren[auto](){\bmu_+ \cdot \d \bn_+[\bV] + \bmu_- \cdot \d \bn_-[\bV]}[\bW]
	\quad
	\text{by \eqref{eq:shape_derivative_log}}
	\notag
	\\
	&
	=
	\d \bmu_+[\bW] \cdot \d \bn_+[\bV]
	+
	\d \bmu_-[\bW] \cdot \d \bn_-[\bV]
	\notag
	\\
	&
	\quad
	+
	\bmu_+ \cdot \d^2 \bn_+[\bV, \bW]
	+
	\bmu_- \cdot \d^2 \bn_-[\bV, \bW]
	,
	\label{eq:d2logmu}
\end{align}
again under the assumption $\d \bV[\bW] = 0$.
Upon close inspection, the symmetry of the above expression in terms of $\bV$ and $\bW$ is not obvious, so we address it in the following lemma.

\begin{lemma}
	\label{lemma:second-order-shape-derivative:establishing-symmetry}
	The following equality holds:
	\begin{multline}
		\label{eq:second-order-shape-derivative:establishing-symmetry}
		\d \bmu_+[\bW] \cdot \d \bn_+[\bV]
		+
		\d \bmu_-[\bW] \cdot \d \bn_-[\bV]
		\\
		=
		\bt_E^\transp \, (D_\Gamma \bW)^\transp \paren[big](){\bmu_+ \bn_+^\transp+ \bmu_- \bn_-^\transp} (D_\Gamma \bV) \bt_E
		.
	\end{multline}
\end{lemma}
\begin{proof}
	We consider two one-sided tangential Jacobians $\paren[auto](){D_\Gamma\bV}_+$ and $\paren[auto](){D_\Gamma\bV}_-$ on the shared edge $E$ between two neighboring triangles $T_+$ and $T_-$.
	That is, $\paren[auto](){D_\Gamma\bV}_+$ stems from extending $D_\Gamma\bV$ from the triangle~$T_+$ to the edge $E$, and likewise for $\paren[auto](){D_\Gamma\bV}_-$.
	Observe that since $\bn_+$, $\bmu_+$ and $\bt_E$ form an orthonormal basis of $\R^3$, we have
	\begin{equation*}
		\id
		=
		\bn_+ \bn_+^\transp
		+
		\bmu_+ \bmu_+^\transp
		+
		\bt_E \bt_E^\transp
		.
	\end{equation*}
	 In view of $\paren[auto](){D_\Gamma \bV}_+ \bn_+ = 0$, we have
	 \begin{align}
		 \d \bn_+[\bV]
		 &
		 =
		 -\paren[big](){ \paren[auto](){D_\Gamma \bV}_+ \id}^\transp \bn_+
		 \quad
		 \text{by \eqref{eq:simple_shape_derivatives:1}}
		 \notag
		 \\
		 &
		 =
		 -\paren[big](){ \paren[auto](){D_\Gamma \bV}_+ \bmu_+ \bmu_+^\transp + \paren[auto](){D_\Gamma \bV}_+ \bt_E \bt_E^\transp}^\transp \bn_+
		 \notag
		 \\
		 &
		 =
		 - \paren[auto](){ \bmu_+ \bmu_+^\transp + \bt_E \bt_E^\transp}^\transp \paren[auto](){D_\Gamma \bV}_+^\transp \bn_+
		 ,
		 \label{eq:reformulated_dnV}
	 \end{align}
	and similarly for $\d \bn_-[\bV]$.
	Taking the shape derivative on both sides of $ \bmu_+ \cdot \bmu_+ = 1$ yields $\bmu_+ \cdot \d \bmu_+[\bW] = 0$.
	Therefore
	\begin{align*}
		\MoveEqLeft
		\d \bmu_+[\bW] \cdot \d \bn_+[\bV]
		\\
		&
		=
		- \d \bmu_+[\bW] \cdot \bt_E \paren[big](){\bt_E^\transp \paren[auto](){D_\Gamma\bV}_+^\transp \bn_+}
		- \d \bmu_+[\bW] \cdot \bmu_+ \paren[big](){\bmu_+^\transp \paren[auto](){D_\Gamma\bV}_+^\transp \bn_+}
		\quad
		\text{by \eqref{eq:reformulated_dnV}}
		\notag
		\\
		&
		=
		- \d \bmu_+[\bW]^\transp \bt_E \paren[big](){\bt_E^\transp \paren[auto](){D_\Gamma \bV}_+^\transp \bn_+}
		.
	\end{align*}
	In the following, we use that $\d \bmu_+[\bW]^\transp \bt_E = \bt_E^\transp\d \bmu_+[\bW]$ and $\bt_E^\transp \paren[auto](){D_\Gamma \bV}_+^\transp \bn_+ = \bn_+^\transp \paren[auto](){D_\Gamma \bV}_+\bt_E$ because both are scalar quantities.
	By plugging in the definition of $\d \bmu_+[\bW]$, equation~\eqref{eq:simple_shape_derivatives:2}, and using $\bmu_+ \cdot \bt_E = 0$ as well as \eqref{eq:simple_shape_derivatives:2}, one obtains
	\begin{align*}
		\MoveEqLeft
		\d \bmu_+[\bW] \cdot \d \bn_+[\bV]
		\\
		&
		=
		-
		\paren[auto](){\bt_E^\transp \d \bmu_+[\bW]} \paren[auto](){ \bn_+^\transp \paren[auto](){D_\Gamma \bV}_+ \bt_E}
		\\
		&
		=
		-
		\paren[Big](){  \bt_E^\transp \paren[auto](){D_\Gamma \bW}_+ \bmu_+- \bt_E^\transp \paren[auto][]{\paren[auto](){D_\Gamma \bW}_+ + \paren[auto](){D_\Gamma \bW}_+^\transp} \bmu_+} \paren[auto](){\bn_+^\transp \paren[auto](){D_\Gamma \bV}_+ \bt_E}
		\\
		&
		=
		\paren[big](){\bt_E^\transp \paren[auto](){D_\Gamma\bW}_+^\transp \bmu_+} \paren[auto](){\bn_+^\transp \paren[auto](){D_\Gamma \bV}_+ \bt_E}
		\\
		&
		=
		\bt_E^\transp \paren[auto](){D_\Gamma \bW}_+^\transp \paren[auto](){\bmu_+ \bn_+^\transp} \paren[auto](){D_\Gamma \bV}_+ \bt_E
		.
		\\
		\intertext{Analogously, one obtains}
		\MoveEqLeft
		\d \bmu_-[\bW] \cdot \d \bn_-[\bV]
		\\
		&
		=
		\bt_E^\transp \paren[auto](){D_\Gamma \bW}_-^\transp (\bmu_- \bn_-^\transp) \paren[auto](){D_\Gamma \bV}_- \bt_E
		.
	\end{align*}
	As before, $(D_\Gamma \bV) \, \bt_E$  is well-defined on $ E$ since $ \paren[auto](){ D_\Gamma \bV}_+ \bt_E = \paren[auto](){ D_\Gamma \bV}_- \bt_E$, and similarly for $ \paren[auto](){D_\Gamma \bW} \bt_E$.
	Therefore, \eqref{eq:second-order-shape-derivative:establishing-symmetry} can be rewritten as
	\begin{equation*}
		\d \bmu_+[\bW] \cdot \d \bn_+[\bV]
		+
		\d \bmu_-[\bW] \cdot \d \bn_-[\bV]
		=
		\bt_E^\transp \paren[auto](){D_\Gamma \bW}^\transp (\bmu_+ \bn_+^\transp+ \bmu_- \bn_-^\transp) \paren[auto](){D_\Gamma \bV} \bt_E.
	\end{equation*}
	The above expression is symmetric, provided that $\bmu_+ \bn_+^\transp + \bmu_- \bn_-^\transp $ is a symmetric matrix, which is not immediately obvious.
	Since $\bn_-$, $\bmu_-$ and $\bt_E$ form an orthonormal basis of $\R^3$, $\id = \bn_- \bn_-^\transp + \bmu_- \bmu_-^\transp + \bt_E \bt_E^\transp $.
	Then, using $\bt_E^\transp \bmu_+ = 0$ and $\bn_+^\transp \bt_E = 0$,
	\begin{align*}
		\bmu_+ \bn_+^\transp + \bmu_- \bn_-^\transp
		&
		=
		\id \, \bmu_+ \bn_+^\transp \, \id
		+
		\bmu_- \bn_-^\transp
		\\
		&
		=
		\paren[big](){\bn_- \bn_-^\transp + \bmu_- \bmu_-^\transp} \, \bmu_+ \bn_+^\transp \paren[big](){\bn_- \bn_-^\transp + \bmu_- \bmu_-^\transp}
		+
		\bmu_- \bn_-^\transp
		\\
		&
		=
		\bn_- \paren[normal](){\bn_-^\transp \bmu_+} \paren[normal](){\bn_+^\transp \bn_-} \, \bn_-^\transp
		+
		\bmu_- \paren[normal](){\bmu_-^\transp \bmu_+} \paren[normal](){\bn_+^\transp \bmu_-} \, \bmu_-^\transp
		\\
		&
		\quad
		+
		\bn_- \paren[normal](){\bn_-^\transp \bmu_+} \paren[normal](){\bn_+^\transp \bmu_-} \, \bmu_-^\transp
		+
		\bmu_- \paren[big][]{(\bmu_-^\transp \bmu_+) (\bn_+^\transp \bn_-) + 1} \, \bn_-^\transp
		.
	\end{align*}
	For the above matrix to be symmetric, it remains to show that
	\begin{equation*}
		(\bn_-^\transp \bmu_+)(\bn_+^\transp \bmu_-)
		-
		(\bmu_-^\transp \bmu_+)(\bn_+^\transp \bn_-)
		=
		1.
	\end{equation*}
	This identity is true by the vector triple product of the cross product since
	\begin{equation*}
		\bn_-^\transp \paren[big][]{\bmu_+ (\bn_+^\transp \bmu_-) - \bn_+ (\bmu_+^\transp \bmu_-)}
		=
		\bn_-^\transp \paren[big][]{\bmu_- \times (\bmu_+ \times \bn_+)}
		=
		\bn_-^\transp \bn_-
		=
		1
		.
	\end{equation*}
\end{proof}

We conclude this subsection by collecting the terms in contributing to the second derivative of the augmented Lagrangian, that is given by
\begin{align}
	\MoveEqLeft
	\d^2 \cL_\rho[\bV, \bW]
	\notag
	\\
	&
	=
	\d^2 \cF[\bV, \bW]
	+
	\sum_{E \in \cE} \d^2\abs{E}_2[\bV,\bW] \, \abs{d_E}
	\notag
	\\
	&
	\quad
	+
	\frac{\rho}{2} \sum_{E \in \cE} \d^2 \abs{E}_2[\bV, \bW] \paren[big][]{d_E - (\logarithm{\bn_+}{\bn_-}) \cdot \bmu_+ - b_E}^2
	\notag
	\\
	&
	\quad
	-
	\frac{\rho}{2} \sum_{E \in \cE} \d^2 \abs{E}_2[\bV, \bW] \, b_E^2
	\notag
	\\
	&
	\quad
	+
	\rho \sum_{E \in \cE} \abs{E}_2 \div_E \bV \, \paren[big][]{d_E - (\logarithm{\bn_+}{\bn_-}) \cdot \bmu_+ - b_E}
	\d \paren[big](){-(\logarithm{\bn_+}{\bn_-}) \cdot \bmu_+}[\bW]
	\notag
	\\
	&
	\quad
	+
	\rho \sum_{E \in \cE} \abs{E}_2 \div_E \bW \, \paren[big][]{d_E - (\logarithm{\bn_+}{\bn_-}) \cdot \bmu_+ - b_E}
	\d \paren[big](){-(\logarithm{\bn_+}{\bn_-}) \cdot \bmu_+}[\bV]
	\notag
	\\
	&
	\quad
	+
	\rho \sum_{E \in \cE} \abs{E}_2\, \paren[big][]{d_E - (\logarithm{\bn_+}{\bn_-}) \cdot \bmu_+ - b_E}
	\d^2 \paren[big](){-(\logarithm{\bn_+}{\bn_-}) \cdot \bmu_+}[\bV, \bW]
	\notag
	\\
	&
	\quad
	+
	\rho \sum_{E \in \cE} \abs{E}_2
	\d \paren[big](){-(\logarithm{\bn_+}{\bn_-}) \cdot \bmu_+}[\bW]
	\d \paren[big](){-(\logarithm{\bn_+}{\bn_-}) \cdot \bmu_+}[\bV]
	.
	\label{eq:d2lagrangian}
\end{align}
For brevity, we leave it to the reader to insert the terms \eqref{eq:shape_derivative_log}, \eqref{eq:d2E}, \eqref{eq:d2logmu} into \eqref{eq:d2lagrangian}.

\subsection{Newton's Method for Shape Optimization}
\label{subsection:Newton_method}

In this section we describe a globalized, inexact Newton method for the shape optimization step (\cref{line:Split_Bregman_discrete_shape_step} of \cref{algorithm:split_Bregman}), \ie, the partial minimization of the augmented Lagrangian~$\cL_\rho$ with respect to the vertex positions of the surface~$\Gamma$.
The inexact Newton direction is obtained from the iterative solution of the shape Newton system
\begin{equation}
	\label{eq:Newton_problem}
	\d^2 \cL_\rho(\Gamma, d, b)[\bW, \bV]
	=
	- \d \cL_\rho(\Gamma, d, b)[\bV]
	\quad
	\text{for all }
	\bV \in \CG{1}(\Gamma, \R^3)
	,
\end{equation}
where the unknown $\bW$ is sought in $\CG{1}(\Gamma, \R^3)$ as well.
This linear system can be represented in the form $A \bw = \bb$ \wrt the standard nodal basis of $\CG{1}(\Gamma, \R^3)$.
It is solved using a preconditioned, truncated conjugate gradient (CG) method (\cref{algorithm:truncated_CG}) inspired by \cite[Algorithm~7.1]{NocedalWright:2006:1}.
The preconditioner $M$ is induced by an incomplete Cholesky decomposition of the matrix respresentation of the~$H^1$-inner product
\begin{equation}
	\label{eq:H1_inner}
	\inner{\bW}{\bV}_\Gamma
	\coloneqq
	\int_\Gamma \bW \cdot \bV \d s + c\int_\Gamma D_\Gamma \bW \dprod D_\Gamma \bV \d s
\end{equation}
on $\CG{1}(\Gamma, \R^3)$ with some constant $c > 0$.
In our experiments we use $c = 10^{-2}$.

\begin{algorithm}
	\caption{Truncated preconditioned CG method for the solution of $A\bw = \bb$.}
	\label{algorithm:truncated_CG}
	\begin{algorithmic}[1]
		\Require matrix $A \in \R^{n \times n}$, symmetric
		\Require matrix $M \in \R^{n \times n}$, symmetric and positive definite
		\Require right-hand side $\bb \in \R^n$
		\Ensure inexact solution to $A\bw = \bb$
		\Ensure $\poscurvature$
		\State Set $k \coloneqq 0$ and $\poscurvature \coloneqq \TRUE$
		\State Set $\bw_0 \coloneqq 0$, $\br_0 \coloneqq -\bb$, $\bp_0 \coloneqq -M^{-1} \br_0$
		\While{$\norm{\br_k}_{M^{-1}} \ge \min(0.5, \sqrt{\norm{\br_0}_{M^{-1}}}) \, \norm{\br_0}_{M^{-1}}$}
		\State Set $\gamma_k \coloneqq \bp_k^\transp A \, \bp_k$
		\If{$\gamma_k \le 0 $}
		\If{$k = 1$}
		\State Set $\poscurvature \coloneqq \FALSE$
		\State \Return $\bp_0$, flag
		\Else
		\State \Return $\bw_k$, flag
		\EndIf
		\EndIf
		\State Set $\alpha_k \coloneqq \norm{\br_k}_{M^{-1}}^2 / \gamma_k$
		\State Set $\bw_{k+1} \coloneqq \bw_k + \alpha_k \, \bp_k$
		\State Set $\br_{k+1} \coloneqq \br_k + \alpha_k \, A \bp_k$
		\State Set $\beta_{k+1} \coloneqq \norm{\br_{k+1}}_{M^{-1}}^2 / \norm{\br_k}_{M^{-1}}^2$
		\State Set $\bp_{k+1} \coloneqq -M^{-1} \br_{k+1} + \beta_{k+1} \, \bp_k$
		\State Set $k \coloneqq k+1$
		\EndWhile
		\State \Return $\bw_k$, flag
	\end{algorithmic}
\end{algorithm}

The truncated CG method starts by computing an approximation of the negative shape gradient direction \wrt the inner product \eqref{eq:H1_inner} and proceeds with CG iterations as long as the search directions remain directions of positive curvature \wrt the shape Hessian.
The latter is guaranteed to always be the case if the shape Hessian is positive definite.
Iterations stop when the desired accuracy is reached, or else when a search direction of non-positive curvature is encountered.
This procedure guarantees that the inexact solution $\bW$ of \eqref{eq:Newton_problem} is a descent direction.
We use the (inexact) Newton search direction returned by \cref{algorithm:truncated_CG} in a descent algorithm globalized by line search (\cref{algorithm:Newton_solve}).
Depending on the termination condition of \cref{algorithm:truncated_CG}, expressed by the flag, we use different initial step sizes for the Armijo line search, treating $\bW$ as either a gradient-type or a Newton-type direction.
Furthermore, we check (see \cref{line:angle_condition} of \cref{algorithm:Newton_solve}) whether the quality of the approximate Newton direction is sufficient compared to the negative gradient direction, as long as the gradient norm is still large.
If not, we fall back to a gradient step.
This is a standard procedure in globalized Newton methods; see, \eg, \cite[p.49]{UlbrichUlbrich:2012:1}.
We choose the $M$-inner product to be the matrix obtained from an incomplete Cholesky decomposition of the matrix respresentation of the~$H^1$-inner product \eqref{eq:H1_inner}.

\begin{algorithm}[htb]
	\caption{Globalized Newton scheme for the approximate minimization of $\cL_\rho(\Gamma, d, b)$ with respect to the vertex position of $\Gamma$.}
	\label{algorithm:Newton_solve}
	\begin{algorithmic}[1]
		\Require initial mesh $\sequence{\Gamma}{0}$
		\Require variations of the normals $d$ and Lagrange multiplier estimate $b$
		\Require Armijo parameters $\sigma \in (0,1)$ and backtracking parameter $\alpha \in (0,1)$
		\Require initial alternate step size $\sequence{\Delta}{0} \in (0,1]$
		\Require maximal number of iterations $\iter_\textup{max}$
		\Require inner product $M$
		\Require $\text{TOL}_\text{smooth}$
		\Ensure optimized mesh $\Gamma$
		\State Set $\ell \coloneqq 0$
		\While{$ \norm{ \nabla_M \cL_\rho}_M \geq \text{TOL}_\text{smooth}$ and $\ell < \iter_\textup{max}$}
		\State Use \cref{algorithm:truncated_CG} to compute an inexact Newton direction $\bW \in \CG{1}(\Gamma, \R^3)$ from \eqref{eq:Newton_problem} (flag is \TRUE), or the negative gradient direction $\bW \in \CG{1}(\Gamma, \R^3)$ (flag is \FALSE)
		\label{line:Newton_solve_compute_newton_dir}
		\If{$\poscurvature$}
			\If{$\d \cL_\rho(\sequence{\Gamma}{\ell}, d, b)[\bW] \leq -\min \{0.1, 10^{-6} \norm{\bW}_M^{0.1}\} \norm{\bW}_M \norm{\nabla_M \cL_\rho}_M$}
			\label{line:angle_condition}
				\State Set $\Delta_\textup{start} \coloneqq 1$
				\Comment{use the inexact Newton direction}
			\Else
			\Comment{use the negative gradient direction}
			\State Set $\Delta_\textup{start} \coloneqq 2 \sequence{\Delta}{\ell}$ and $\bW \coloneqq - \nabla_M \cL_\rho(\sequence{\Gamma}{\ell}, d, b)$
			\EndIf
		\Else
		\Comment{use the negative gradient direction}
			\State Set $\Delta_\textup{start} \coloneqq 2 \sequence{\Delta}{\ell}$
		\EndIf
		\State Find $\sequence{\Delta}{\ell+1}$ by backtracking, starting at $\Delta_\textup{start}$, such that $\cL_\rho(\sequence{\Gamma}{\ell}_{\sequence{\Delta}{\ell+1}}[\bW], d, b) \le \cL_\rho(\sequence{\Gamma}{\ell}, d, b) + \sigma \, \sequence{\Delta}{\ell+1} \d\cL_\rho(\sequence{\Gamma}{\ell}, d, b)[\bW]$ holds
		\label{line:Newton_solve_Armijo_gradient}
		\State Set $\sequence{\Gamma}{\ell+1} \coloneqq \sequence{\Gamma}{\ell}_{\sequence{\Delta}{\ell+1}}[\bW] \coloneqq \setDef[auto]{x + \sequence{\Delta}{\ell+1} \, \bW(x)}{x \in \sequence{\Gamma}{\ell}}$
		\State Set $\ell \coloneqq \ell+1$
		\EndWhile
	\end{algorithmic}
\end{algorithm}

Let us emphasize that in later iterations of the split Bregman iteration (\cref{algorithm:split_Bregman}), one can expect to use exclusively Newton directions and converge superlinearly to a local minimizer of the smooth subproblem.
In practice this is indeed observed after only a few iterations of \cref{algorithm:split_Bregman}.

As initial step size, we use $\Delta_\textup{start} = 1$ in case of a Newton-type direction, or twice the most recent successful value of $\Delta_\textup{start}$ from the previous call \eqref{eq:Newton_problem} in case of a gradient-type direction.
To prevent breaking the mesh by self-intersections or degenerate cells in the first few iterations, we are using a rather small initial step size $\Delta_\textup{start}$ when \cref{algorithm:Newton_solve} is called in the very first iteration of \cref{algorithm:split_Bregman}.

\subsection{Convergence Critera for the Split Bregman Method}
\label{subsection:convergence_criteria}

In this subsection we discuss a convergence criterion for \cref{algorithm:split_Bregman}, based on residual norms.
Due to the shape dependencies of the constraint penalties, our setting is slightly more complex than \cite[Section~3.2]{BoydParikhChuPeleatoEckstein:2010:1}.
The residuals are derived from the optimality conditions of the (non-augmented) Lagrangian
\begin{equation}
	\label{eq:Lagrangian_original}
	\cL(\Gamma, d, b)
	\coloneqq
	\cF(\Gamma)
	+
	\beta \sum_{E \in \cE} \abs{E}_2 \, \abs{d_E}
	+
	\rho \sum_{E \in \cE} \abs{E}_2 \, b_E \paren[big][]{(\logarithm{\bn_+}{\bn_-}) \cdot \bmu_+ - d_E}
	,
\end{equation}
where $b$ is the scaled Lagrange multiplier, that is connected to the true Lagrange multiplier $\lambda$ by $\lambda_E = \rho \, b_E$.
The first-order optimality conditions for an optimal point $(\Gamma^*, d^*, b^*)$ are then
\begin{align*}
	0
	&
	= \d \cL(\Gamma^*, d^*, b^*)[\bV]
	\quad
	\text{for all }
	\bV \in \CG{1}(\Gamma, \R^3)
	,
	\\
	0
	&
	\in \partial_d \cL (\Gamma^*, d^*, b^*)
	,
	\\
	0
	&
	=
	\partial_b \cL(\Gamma^*, d^*, b^*)
	,
\end{align*}
where $\partial_d$ denotes the convex sub-derivative \wrt $d$ and $\partial_b$ denotes the derivative \wrt $b$.
This motivates us to monitor the following residuals
\begin{align*}
	\sequence{q}{k+1}[\bV]
	&
	\coloneqq
	\d \cL(\sequence{\Gamma}{k+1}, \sequence{d}{k+1}, \sequence{b}{k+1})[\bV]
	,
	\\
	\sequence{r}{k+1}
	&
	\coloneqq
	\partial_b \cL(\sequence{\Gamma}{k+1}, \sequence{d}{k+1}, \sequence{b}{k+1})
	.
\end{align*}
We derive a suitable third residual $\sequence{s}{k+1} \in \partial_d \cL(\sequence{\Gamma}{k+1}, \sequence{d}{k+1}, \sequence{b}{k+1})$, an element of the convex subdifferential, analogously to \cite{BoydParikhChuPeleatoEckstein:2010:1}, \cref{line:Split_Bregman_discrete_tv_step} in \cref{algorithm:split_Bregman} ensures $\bnull\in\partial_d \cL_\rho(\sequence{\Gamma}{k}, \sequence{d}{k+1}, \sequence{b}{k})$, \ie,
\begin{align*}
	0
  &\in \beta \sum_E \abs{E}_2 \partial \abs{\sequence{d_E}{k+1}}
	+
  \rho \sum_{E \in \cE}
  \paren[auto](){
	\sequence{d_E}{k+1}
	-
	(\logarithm{\sequence{\bn_+}{k}}{\sequence{\bn_-}{k}}) \cdot \sequence{\bmu_+}{k}
	-
	\sequence{b_E}{k}
  }
  .
\end{align*}
Plugging in the multiplier update (\cref{line:Split_Bregman_discrete_multiplier_step}) and inserting \eqref{eq:Lagrangian_original} yields
\begin{equation*}
  \bnull
  \in
  \partial_d \cL(\sequence{\Gamma}{k+1},\sequence{d}{k+1},\sequence{b}{k+1})
  - \sequence{s}{k+1}
  ,
\end{equation*}
where
\begin{equation*}
  \sequence{s}{k+1}
  \coloneqq
  -\rho
  \sum_{E \in \cE}
  (\logarithm{\sequence{\bn_+}{k+1}}{\sequence{\bn_-}{k+1}}) \cdot \sequence{\bmu_+}{k+1}
	-
	(\logarithm{\sequence{\bn_+}{k}}{\sequence{\bn_-}{k}}) \cdot \sequence{\bmu_+}{k}
  .
\end{equation*}
The residuals $\sequence{r}{k+1}$ and $\sequence{s}{k+1}$ are known from the literature as primal and dual residuals, respectively; see, \eg, \cite{BoydParikhChuPeleatoEckstein:2010:1} for details.
On an edge~$E$, they are evaluated according to
\begin{subequations}
	\begin{align}
		\label{eq:dual_residual_edge}
		\sequence{s_E}{k+1}
		&
		\coloneqq
		- \rho \paren[auto](){(\logarithm{\sequence{\bn_+}{k+1}}{\sequence{\bn_-}{k+1}}) \cdot \sequence{\bmu_+}{k+1} - (\logarithm{\sequence{\bn_+}{k}}{\sequence{\bn_-}{k}}) \cdot \sequence{\bmu_+}{k}}
		,
		\\
		\label{eq:primal_residual_edge}
		\sequence{r_E}{k+1}
		&
		\coloneqq
		\paren[big](){\logarithm{\sequence{\bn_+}{k+1}}{\sequence{\bn_-}{k+1}}} \cdot \sequence{\bmu_+}{k+1} - \sequence{d_E}{k+1}
		.
	\end{align}
\end{subequations}
Overall we will consider the residual norms
\begin{subequations}
	\label{eq:norms_residuals}
	\begin{align}
		\norm{\sequence{r}{k+1}}_{L^2}^2
		&
		\coloneqq
		\sum_{E \in \cE} \abs{E}_2 \abs{\sequence{r_E}{k+1}}^2
		\quad
		\text{from \eqref{eq:primal_residual_edge}}
		,
		\label{eq:alg1_primalres}
		\\
		\norm{\sequence{s}{k+1}}_{L^2}^2
		&
		\coloneqq
		\sum_{E \in \cE} \abs{E}_2  \abs{\sequence{s_E}{k+1}}^2
		\quad
		\text{from \eqref{eq:dual_residual_edge}}
		,
		\label{eq:alg1_dualres}
		\\
		\norm{\sequence{q}{k+1}}_{(H^1)^*}
		&
		\coloneqq
		\sup \setDef[big]{\sequence{q}{k+1}[\bV]}{\bV \in \CG{1}(\sequence{\Gamma}{k+1}, \R^3), \; \norm{\bV}_{H^1} \le 1}
		.
		\label{eq:alg1_shaperes}
	\end{align}
\end{subequations}

In common ADMM frameworks, the third residual $\sequence{q}{k+1}$, \ie, the derivative \wrt to the primal variable that was updated most recently (in our case the vertex positions describing~$\Gamma$), vanishes automatically by design.
This is, however, not the case in the shape optimization problem here, because the penalty term on an edge is scaled by the respective edge length.
In order to derive a more convenient representation of the third residual $\sequence{q}{k+1}[\bV]$, we utilize that in case of an exact minimization in \cref{line:Split_Bregman_discrete_shape_step} of \cref{algorithm:split_Bregman}, the derivative $\d\cL_\rho(\Gamma^{k+1}, d^{k+1}, b^k) = 0$, see \eqref{derivative_augmented_Lagrangian}.
Plugging in the multiplier update (\cref{line:Split_Bregman_discrete_multiplier_step} of \cref{algorithm:split_Bregman}) into this condition yields
	\begin{align*}
	0
	&
	=
	\d \cF(\sequence{\Gamma}{k+1})[\bV]
	+
	\beta \sum_{E \in \cE} \abs{E}_2 \div_E \bV \abs{\sequence{d_E}{k+1}}
	\\
	&
	\quad
	-
	\frac{\rho}{2} \sum_{E \in \cE} \abs{E}_2 \div_E \bV \paren[big][]{(\logarithm{\bn_+}{\bn_-}) \cdot \bmu_+ - \sequence{d_E}{k+1}}^2
	\\
	&
	\quad
	+
	\rho \sum_{E \in \cE} \abs{E}_2 \div_E \bV \, \sequence{b_E}{k+1} \paren[big][]{(\logarithm{\bn_+}{\bn_-}) \cdot \bmu_+ - \sequence{d_E}{k+1}}
	\\
	&
	\quad
	+
	\rho \sum_{E \in \cE} \abs{E}_2 (-\sequence{b_E}{k+1}) \d \paren[auto](){-(\logarithm{\bn_+}{\bn_-}) \cdot \bmu_+}[\bV]
	.
\end{align*}
Plugging in the definition of $\d \cL(\sequence{\Gamma}{k+1}, \sequence{d}{k+1}, \sequence{b}{k+1})[\bV]$, we have
\begin{equation*}
	\d \cL(\sequence{\Gamma}{k+1}, \sequence{d}{k+1}, \sequence{b}{k+1})[\bV]
	=
	\frac{\rho}{2} \sum_{E \in \cE} \abs{E}_2 \div_E \bV \paren[big][]{(\logarithm{\bn_+}{\bn_-}) \cdot \bmu_+ - \sequence{d_E}{k+1}}^2
\end{equation*}
and thus
\begin{equation}
	\sequence{q}{k+1}[\bV]
	=
	\frac{\rho}{2} \sum_{E \in \cE} \abs{E}_2 \div_E \bV \paren[big][]{(\logarithm{\bn_+}{\bn_-}) \cdot \bmu_+ - \sequence{d_E}{k+1}}^2
	.
\end{equation}
Consequently, $\sequence{q}{k+1}$ vanishes if $\sequence{r}{k+1}$ does, provided that the shape subproblem is solved sufficiently accurately by the Newton scheme.
Therefore, we only monitor $\sequence{r}{k+1}$ and $\sequence{s}{k+1}$ to adapt the penalty parameter $\rho$ in the next subsection.

\subsection{Adaptive Penalty Parameter Selection}
\label{subsection:penalty_parameter_selection}

Based on the primal and dual residuals, it is possible to adapt the penalty parameter~$\rho$ to accelerate \cref{algorithm:split_Bregman}.
Here we follow a commonly used approach described in \cite[Section~3.4.1]{BoydParikhChuPeleatoEckstein:2010:1}.
It utilizes the fact that the dual residual~$s$ is proportional to the penalty parameter $\rho$, see \eqref{eq:dual_residual_edge}.
Therefore, by decreasing~$\rho$, one can reduce the dual residual.
Furthermore, since $\rho$ penalizes the constraint, the primal residual of the next iteration can be reduced by increasing~$\rho$.
This behavior can be utilized to balance both residuals to the same order of magnitude, which yields faster overall convergence.

The strategy is overall described by
\begin{equation*}
	\rho_\textup{new}
	=
	\begin{cases}
		1.5\rho
		&
		\text{if }
		\norm{\sequence{r}{k+1}}_{L^2}
		>
		5 \, \norm{\sequence{s}{k+1}}_{L^2}
		\\
		\rho/1.5
		&
		\text{if }
		\norm{\sequence{r}{k+1}}_{L^2}
		<
		5 \, \norm{\sequence{s}{k+1}}_{L^2}
		\\
		\rho
		&
		\text{else}
		.
	\end{cases}
\end{equation*}
Notice that the scaled Lagrange multipliers $b = \lambda/\rho$ have to be rescaled such that the underlying multiplier estimate $\sequence{\lambda}{k}$ remains unchanged.
As is emphasized in  \cite{BoydParikhChuPeleatoEckstein:2010:1}, the major advantage of this is that a manual selection of the penalty parameter is no longer necessary.

\section{Mesh Denoising Problem}
\label{section:Mesh_Denoising_Problem}

In this section, we consider a mesh denoising problem using the total variation of the normal \eqref{eq:Discrete_TV_log_novel} as a regularizer.
As in most applications of variational mesh denoising that do not use normal filtering, we use an $\ell^2$ fidelity term
\begin{equation}
	\label{eq:mesh_denoising_l1_fidelity}
	\frac{1}{2} \sum_{\vertex \in \cV} \abs{x_\vertex - x^\text{data}_\vertex}_2^2
	.
\end{equation}
Each summand measures the squared distance between the spatial coordinate $x_\vertex$ of vertex~$\vertex$ and the spatial coordinate $x^\text{data}_\vertex$ of the respective vertex on the noisy mesh.

Simply combining this fidelity term with the regularizer can, however, lead to problems because mesh quality is a concern of neither.
Typically, degenerate meshes occur when triangles become too small.
This can be prevented by adding an additional barrier term.
We use here a simple approach and propose
\begin{equation*}
	\cF(\Gamma)
	=
	\sum_{\vertex \in \cV} \abs{x_\vertex - x^\text{data}_\vertex}_2^2
	+
	\tau \sum_{T\in \cT} \frac{1}{\abs{T}}
	.
\end{equation*}
The idea of adding a term for the purpose of maintaining mesh quality and ensuring the existence of discrete solution is not new.
The authors of \cite{WuZhengCaiFu:2015:1}, who used an approximation of the regularizer \eqref{eq:Discrete_TV_log_novel}, employed a $\emph{fairness}$ term in order to prevent mesh folding and related effects.
A more involved barrier term for planar meshes was studied in \cite{HerzogLoayzaRomero:2023:1}.

Overall, we consider here the following problem
\begin{equation}
	\label{eq:mesh_denoising_problem}
	\text{Minimize}
	\quad
	\frac{1}{2} \sum_{\vertex \in \cV} \abs{x_\vertex - x^\text{data}_\vertex}_2^2
	+
	\tau \sum_{T \in \cT} \frac{1}{\abs{T}}
	+
	\beta \sum_{E \in \cE} \abs{E}_2 \, \abs[big]{(\logarithm{\bn_+}{\bn_-}) \cdot \bmu_+}
	,
\end{equation}
where $\beta$ and $\tau$ are two regularization parameters.
The values of $\beta$ and $\tau$ are to be balanced in a way so that the total variation regularizer term dominates unless triangles become extremely small.
\begin{figure}[htb]
	\centering
	\includegraphics[width = 0.32\linewidth]{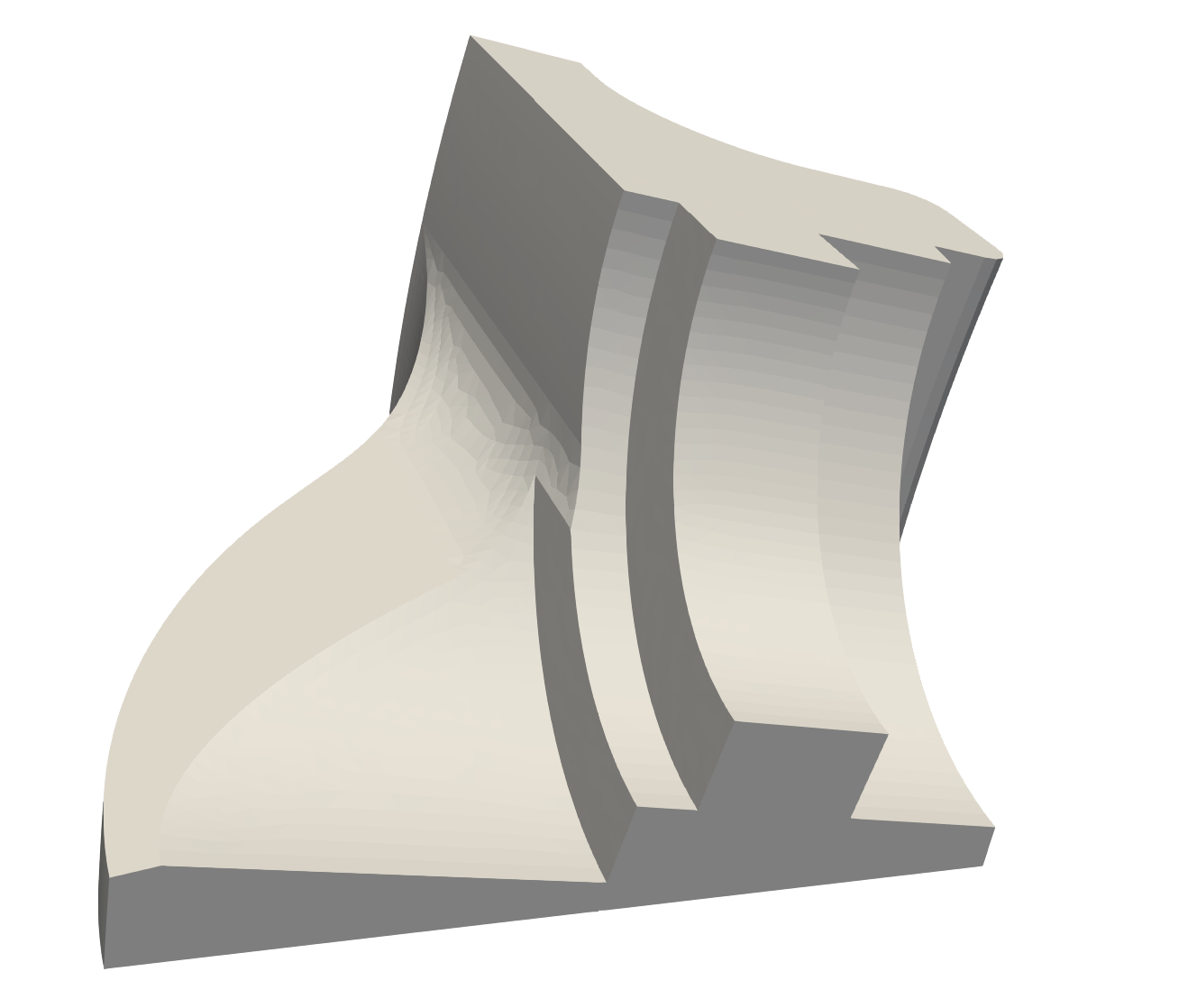}
	\hfill
	\includegraphics[width = 0.32\linewidth]{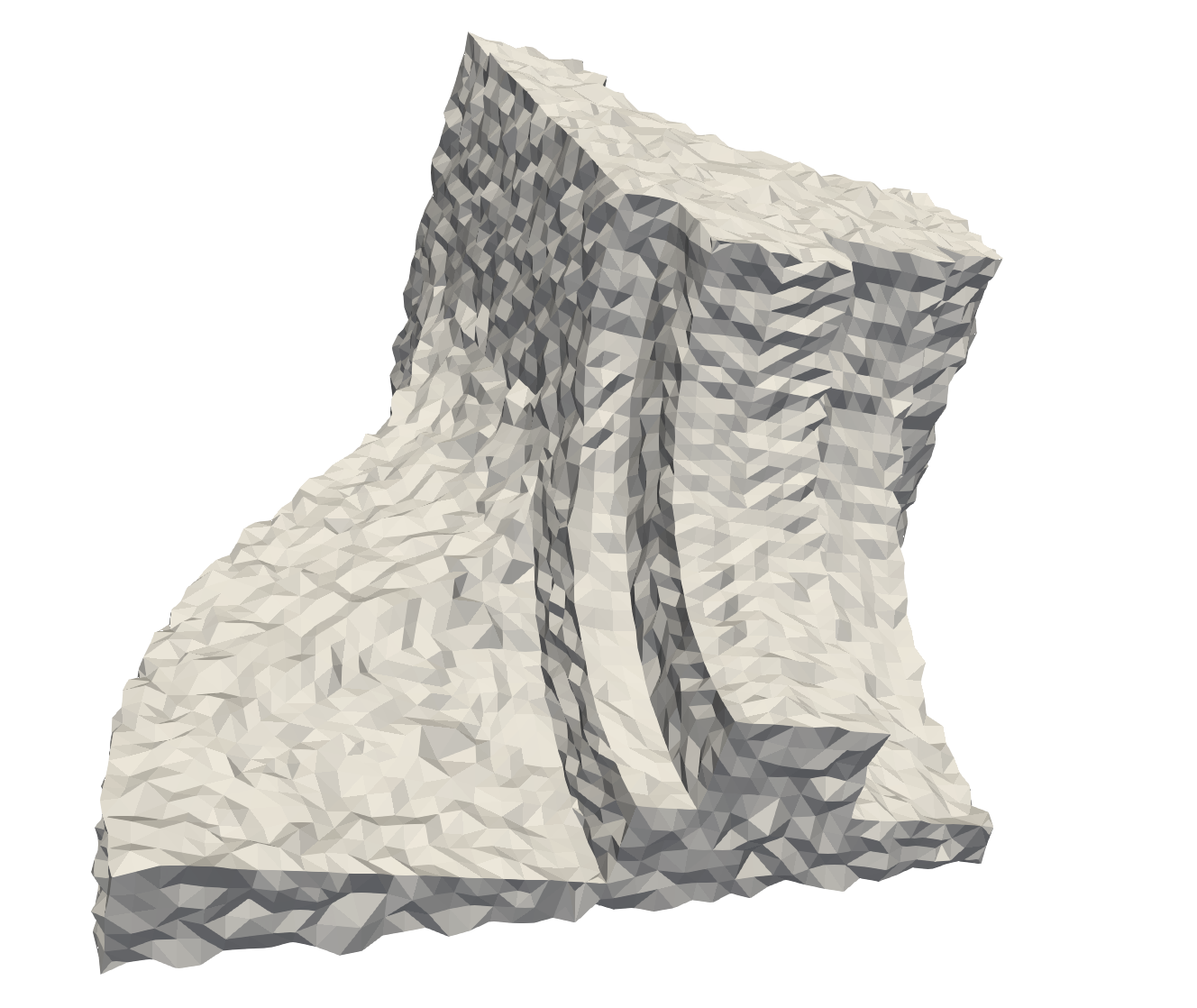}
	\hfill
	\includegraphics[width = 0.32\linewidth]{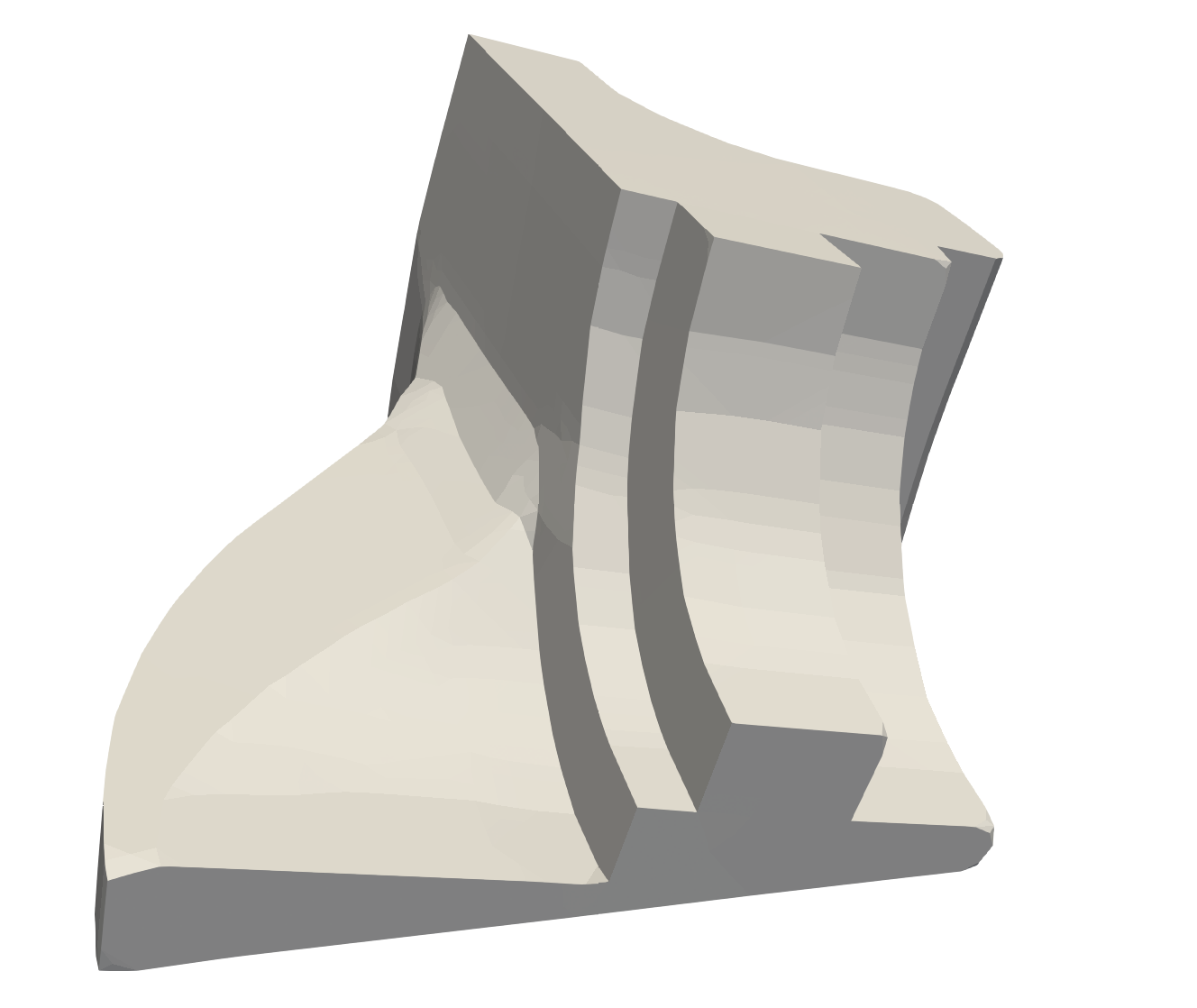}
	\caption{Mesh denoising result using \cref{algorithm:split_Bregman} applied to problem \eqref{eq:mesh_denoising_problem} with $\beta = 2\cdot 10^{-2}$, $\tau = 10^{-8}$ and initial penalty parameter $\rho = 10^{-2}$. Original geometry (left), noisy geometry (middle) and reconstruction (right).}
	\label{figure:Mesh_denoising_using_TV}
\end{figure}

\Cref{figure:Mesh_denoising_using_TV} shows a denoising result obtained using the model \eqref{eq:mesh_denoising_problem} and \cref{algorithm:split_Bregman} after \num{1000}~iterations.
The data $x^\text{data}_\vertex$ is taken from the well-known \emph{fandisk} benchmark problem~\cite{HoppeDeRoseDuchampHalsteadJinMcDonaldSchweitzerStuetzle:1994:1} consisting of $\num{12954}$ triangles, available at the Wolfram Data Respository~\cite{Shedelbower:2022:1}.

Noise is added at each vertex by displacing each vertex in normal direction using a Gaussian distribution with standard deviation $0.2$ times the average length of the edges adjacent to the vertex.
In this and all following experiments we use typical Armijo parameters of $\alpha = 0.5$, $\sigma = 0.5$ as well as $\sequence{\Delta}{0} = 10^{-4}$.
Furthermore, $\iter_\textup{max} = 3$ and $\text{TOL}_\text{smooth}=10^{-3}$ are used in \cref{algorithm:Newton_solve}.
Except for the first iterations, a smaller number of Newton steps are actually required to achieve convergence of the smooth subproblem to the desired tolerance.
The number of Newton steps is visualized in \cref{figure:inner_iterations_rho}, that also contains a plot of the evolution of the adaptive penalty parameter~$\rho$.
Our implementation was done in \fenics (version~$2019.2.\textup{dev}0$).

\begin{figure}
	\centering
	\begin{scaletikzpicturetowidth}{0.45\textwidth}
		\begin{tikzpicture}[baseline, scale = \tikzscale]
			\begin{axis}[
				xlabel = {iteration},
				ylabel = {number of inner iterations},
				grid = major,
				]
				\addplot[thick, each nth point = 2, color = blue] table[x index = 8,y index = 7] {src/Fandisk/Newton/history.txt};
			\end{axis}
		\end{tikzpicture}
	\end{scaletikzpicturetowidth}
	\hfill
	\begin{scaletikzpicturetowidth}{0.45\textwidth}
		\begin{tikzpicture}[baseline, scale = \tikzscale]
			\begin{axis}[
				xlabel = {iteration},
				ylabel = {penalty parameter $\rho$},
				grid = major,
				]
				\addplot[thick, each nth point = 2, color = blue] table[x index = 8,y index = 4] {src/Fandisk/Newton/history.txt};
			\end{axis}
		\end{tikzpicture}
	\end{scaletikzpicturetowidth}
	\caption{Visualization of the number of iterations used for the globalized Newton scheme (\cref{algorithm:Newton_solve}) in each iteration of \cref{algorithm:split_Bregman} (left) applied to the \emph{fandisk} denoising problem described in \cref{section:Mesh_Denoising_Problem}. Evolution of the penalty parameter~$\rho$ over the iterations (right).}
	\label{figure:inner_iterations_rho}
\end{figure}

We use this example to compare our \cref{algorithm:split_Bregman} to the Riemannian split Bregman algorithm \cite{BergmannHerrmannHerzogSchmidtVidalNunez:2020:2}.
Recall that the latter did not use Newton steps for the shape optimization subproblem, and it was based on the total variation formulation \eqref{eq:Discrete_TV_log}, that required parallel transport of the Lagrange multiplier estimates.
Also, since the shape update in \cite{BergmannHerrmannHerzogSchmidtVidalNunez:2020:2} is performed prior to the update of the $\widehat \bd$~variable, the iteration is \emph{dual-feasible} and the residual \wrt $\widehat \bd$ vanishes by design.
Therefore, only the primal and the dual residual are necessary to measure the state of convergence.
The former measures the feasibility, see \eqref{eq:contraint_riemannian_admm}, and the latter measures stationary \wrt shape variations.
We consider
\begin{subequations}
	\label{eq:norms_residuals_RM}
	\begin{align}
		\abs[big]{\sequence{r_\textup{RM}}{k+1}}_{L^2}^2
		&
		\coloneqq
		\sum_{E \in \cE} \abs{E}_2 \abs[big]{\sequence{\widehat \bd_E}{k+1} - \logarithm{\sequence{\bn_+}{k+1}}{\sequence{\bn_-}{k+1}}}_2^2
		,
		\label{eq:algadmm_primalres}
		\\
		\abs[big]{\sequence{s_\textup{RM}}{k+1}}_{(H^1)^*}
		&
		\coloneqq
		\sup \setDef[auto]{\d \cL_\textup{RM}(\sequence{\Gamma}{k+1}, \sequence{\widehat \bd}{k+1}, \sequence{\widehat \bb}{k+1})[\bV]}{%
			\begin{aligned}
				&
				\bV \in \CG{1}(\sequence{\Gamma}{k+1}, \R^3)
				,
				\\
				&
				\norm{\bV}_{H^1} \le 1
			\end{aligned}
		}
		\label{eq:algadmm_dualres}
	\end{align}
\end{subequations}
as a basis of the convergence criterion for the Riemannian ADMM from \cite{BergmannHerrmannHerzogSchmidtVidalNunez:2020:2}.

We furthermore seek to emphasize the effect of Newton's method for the smooth subproblem by comparing \cref{algorithm:split_Bregman} to a variant using only shape gradient directions for the mesh update~$\bW$ in \cref{algorithm:Newton_solve}.
In both cases, the Riemannian ADMM from \cite{BergmannHerrmannHerzogSchmidtVidalNunez:2020:2} as well as the modified version of \cref{algorithm:split_Bregman}, we perform $3$~gradient steps in the shape optimization step using the $H^1$-inner product \eqref{eq:H1_inner} to obtain the gradient direction.
An Armijo rule with the same parameters as before is used; see \cref{line:Newton_solve_Armijo_gradient} of \cref{algorithm:Newton_solve}.

\begin{figure}
	\centering
	\pgfplotsset{
		legend style = {
			at = {(0.01,0.01)},
			anchor = south west,
		},
	}
	\begin{scaletikzpicturetowidth}{0.45\textwidth}
		\begin{tikzpicture}[baseline, scale = \tikzscale]
			\begin{semilogyaxis}[
				xlabel = {iterations},
				ylabel = {combined residual norms},
				grid = major,
				legend entries = {\cref{algorithm:split_Bregman} with Newton, \cref{algorithm:split_Bregman} with $H^1$-gradient, Riemannian ADMM},
				legend style = {font = \tiny},
				legend cell align = {left},
				legend pos = north east,
				],
				\addplot[thick, mark = none, each nth point = 1, color = blue] table[x index = 8, y index = 3] {src/Fandisk/Newton/history.txt};
				\addplot[thick, mark = none, each nth point = 1, color = brown] table[x index = 8,y index = 3] {src/Fandisk/GradientDescent/history.txt};
				\addplot[thick, mark = none, each nth point = 1, color = red] table[x index = 8, y index = 3]  {src/Fandisk/Riemanian/history.txt};
			\end{semilogyaxis}
		\end{tikzpicture}
	\end{scaletikzpicturetowidth}
	\hfill
	\begin{scaletikzpicturetowidth}{0.45\textwidth}
		\begin{tikzpicture}[baseline, scale = \tikzscale]
			\begin{semilogyaxis}[
				xlabel = {CPU time},
				ylabel = {combined residual norms},
				grid = major,
				legend entries = {\cref{algorithm:split_Bregman} with Newton, \cref{algorithm:split_Bregman} with $H^1$-gradient, Riemannian ADMM},
				legend style = {font = \tiny},
				legend cell align = {left},
				legend pos = north east,
				],
				\addplot[thick, mark = none, color = blue] table[x index = 5, y index = 3] {src/Fandisk/Newton/history.txt};
				\addplot[thick, mark = none, color = brown] table[x index =5, y index = 3] {src/Fandisk/GradientDescent/history.txt};
				\addplot[thick, mark = none, color = red] table[x index = 5, y index = 3]  {src/Fandisk/Riemanian/history.txt};
			\end{semilogyaxis}
		\end{tikzpicture}
	\end{scaletikzpicturetowidth}
	\caption{Combined residual norms \eqref{eq:residuals_alg}, \eqref{eq:residuals_RM} for the \emph{fandisk} denoising problem over iterations (left) as well as combined residual norms over CPU time on an AMD Ryzen~5 3600 desktop CPU (right). For all three methods, parameters $\beta = 2\cdot 10^{-2}$ and $\tau = 10^{-8}$ are chosen. For \cref{algorithm:split_Bregman}, the adaptive penalty parameter selection is used starting with $\rho = 10^{-2}$. The Riemannian ADMM method uses the constant penalty parameter $\rho = 1$ for all iterations. Every $20$th iteration is plotted.}
	\label{figure:residuals}
\end{figure}

\Cref{figure:residuals} shows combined norm of the residuals for both variants of \cref{algorithm:split_Bregman} as well as the Riemannian ADMM from \cite{BergmannHerrmannHerzogSchmidtVidalNunez:2020:2}.
For the former, we show the square root of
\begin{equation}
	\label{eq:residuals_alg}
	\sequence{e}{k+1}
	\coloneqq
	\norm[big]{\sequence{r}{k+1}}_{L^2}^2+ \norm[big]{\sequence{s}{k+1}}_{L^2}^2 + \norm[big]{\sequence{q}{k+1}}_{(H^1)^*}^2
\end{equation}
with the summands coming from \eqref{eq:norms_residuals}.
For the Riemannian ADMM, we show the square root of
\begin{equation}
	\label{eq:residuals_RM}
	\sequence{e_\textup{RM}}{k+1}
	\coloneqq
	\norm[big]{\sequence{r_\textup{RM}}{k+1}}_{L^2}^2 + \norm[big]{\sequence{s_\textup{RM}}{k+1}}_{(H^1)^*}^2
\end{equation}
with the summands defined in \eqref{eq:norms_residuals_RM}.
Clearly, \cref{algorithm:split_Bregman} with Newton outperforms both other algorithms.
Indeed, the slow convergence of the two other algorithms can be attributed to the slow reduction of the residuals associated with the respective shape optimization subproblem.

\begin{figure}[htp]
	\centering
	\includegraphics[width = 0.32\linewidth, trim = {15cm 0cm 10cm 0cm}, clip]{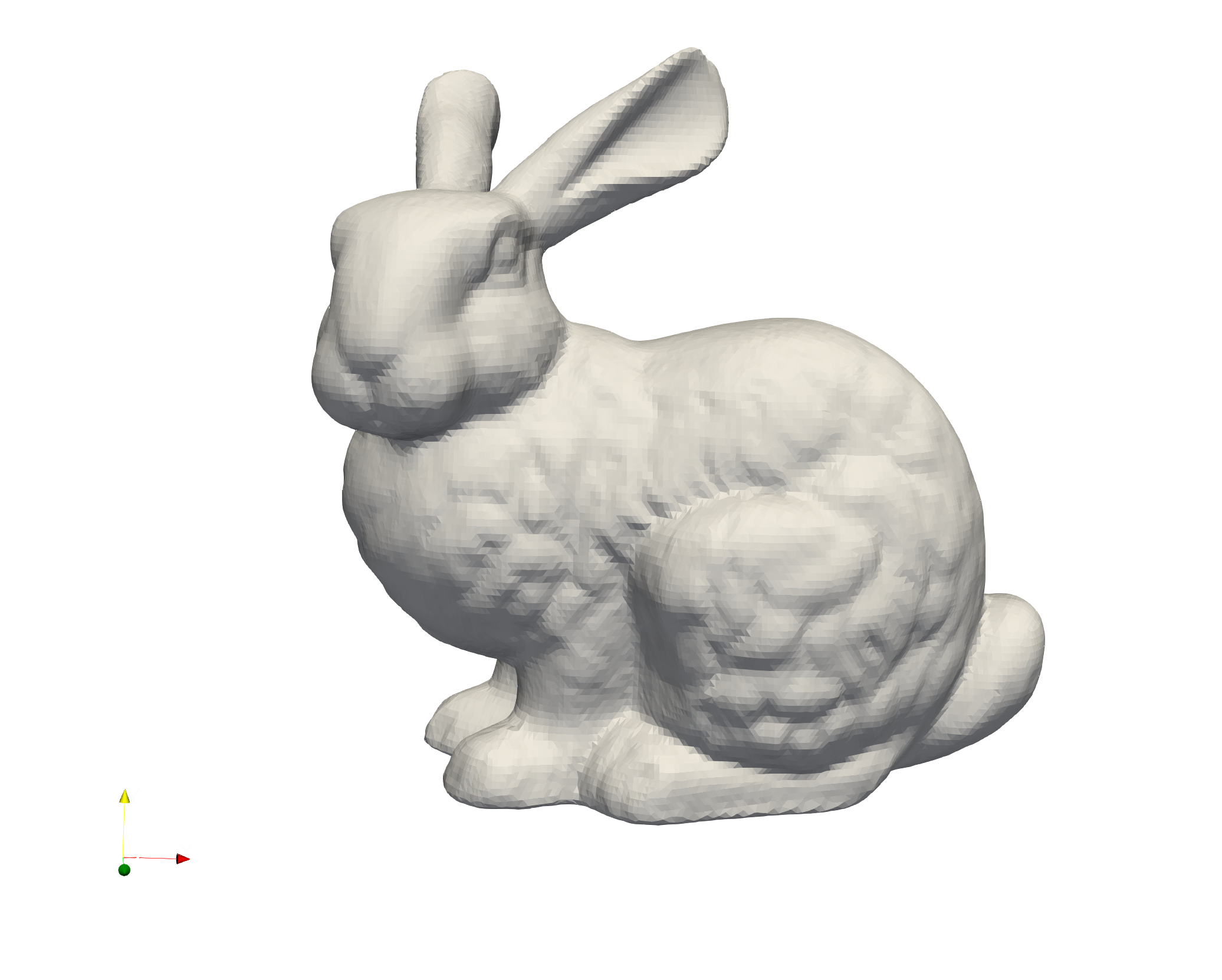}
	\hfill
	\includegraphics[width = 0.32\linewidth, trim = {15cm 0cm 10cm 0cm}, clip]{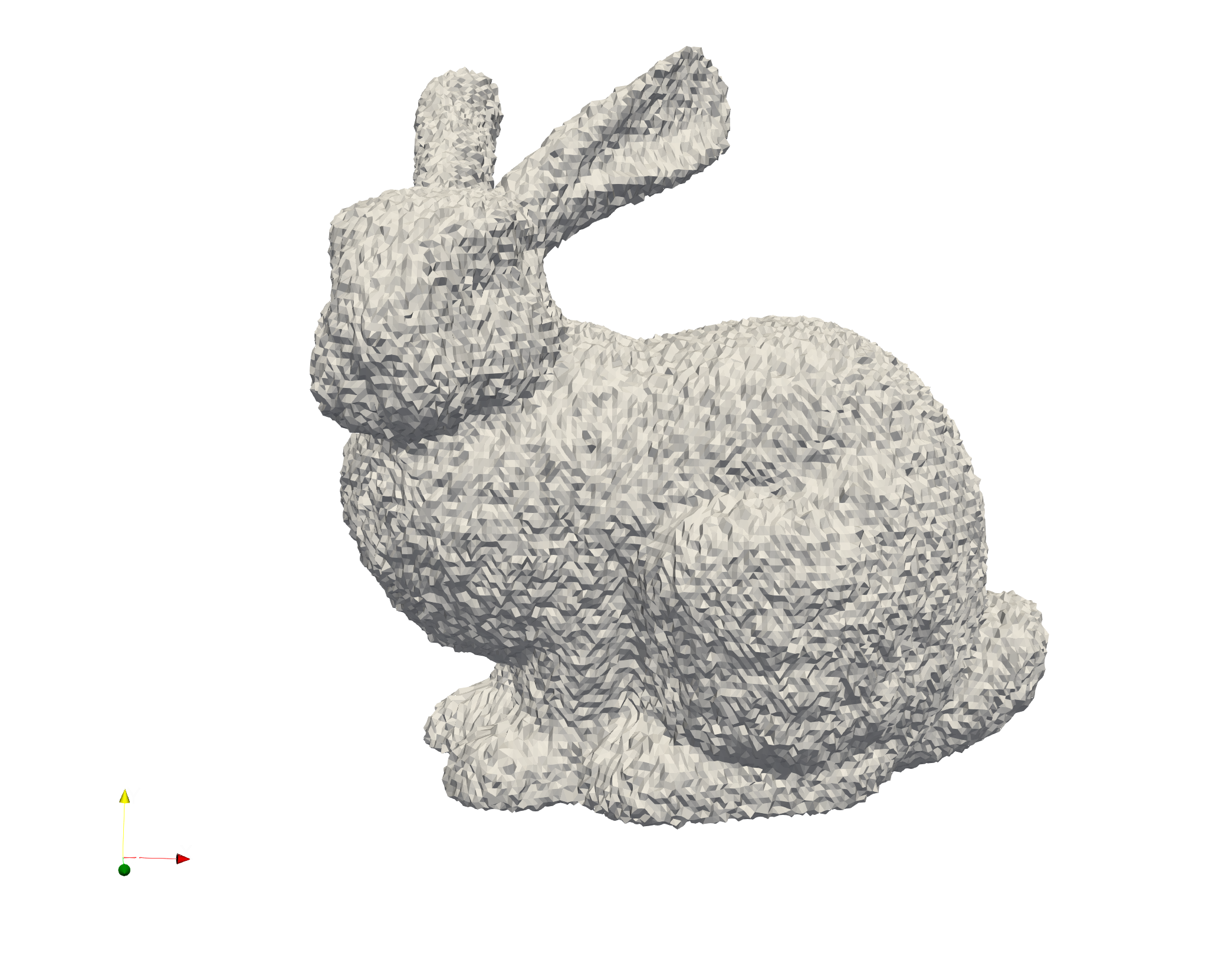}
	\hfill
	\includegraphics[width = 0.32\linewidth, trim = {15cm 0cm 10cm 0cm}, clip]{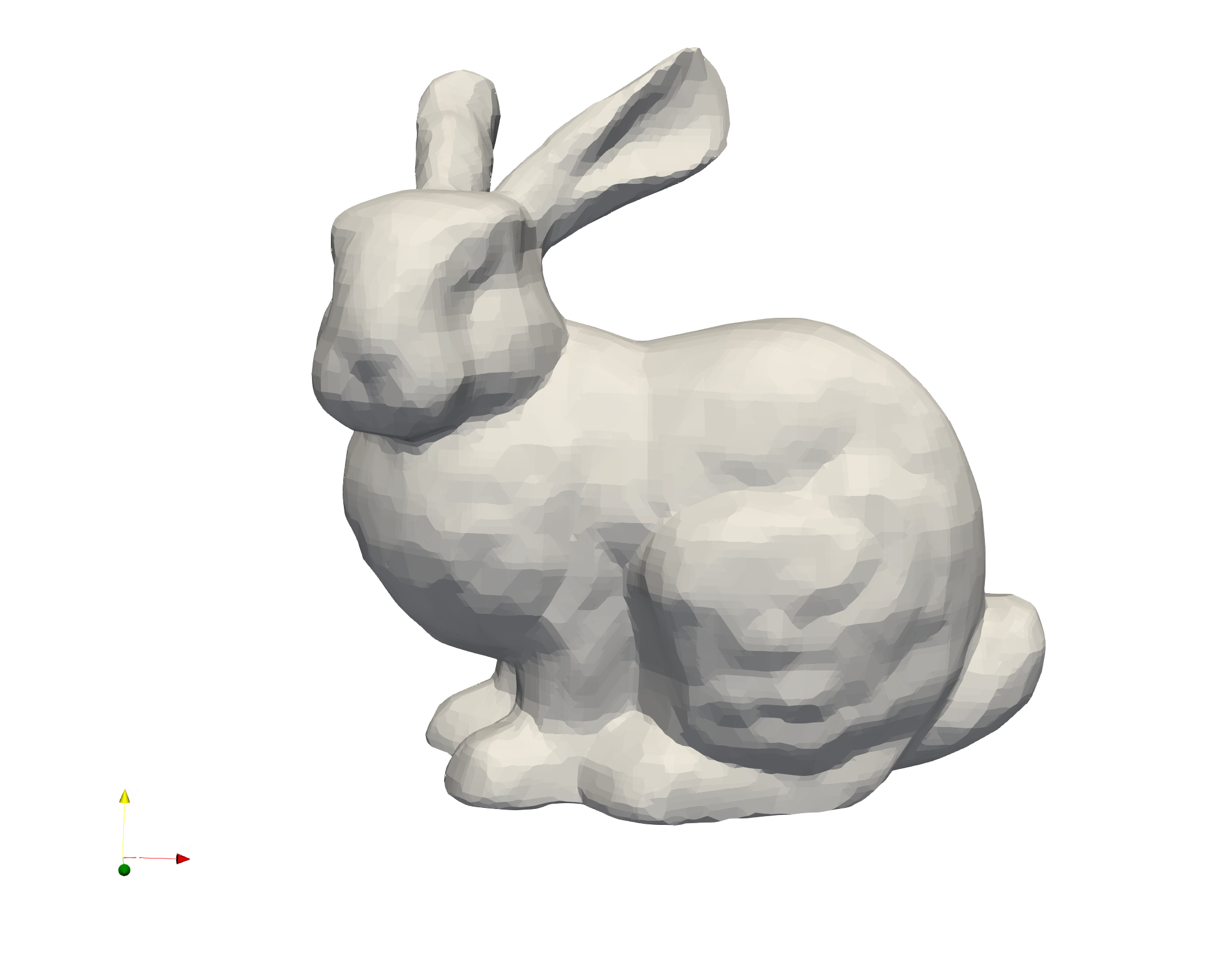}
	\caption{Mesh denoising results using \cref{algorithm:split_Bregman} on problem \eqref{eq:mesh_denoising_problem} with $\beta = 5 \cdot 10^{-3}$, $\tau = 10^{-8}$ and initial penalty parameter $\rho = 10^{-3}$. Original geometry (left), noisy geometry (middle) and reconstruction (right).}
	\label{figure:Mesh_denoising_using_TV_Bunny}
\end{figure}

The denoising results on the \emph{fandisk} geometry in \cref{figure:Mesh_denoising_using_TV} show excellent results in preserving kinks.
However, as expected for TV-based regularizations, there is a mild staircasing in regions of constant curvature.
To this end, we also solve the denoising problem \eqref{eq:mesh_denoising_problem} on the \emph{Stanford bunny} mesh, a geometry with a significantly more intricate structure, that is expected to be more susceptible to this deficiency.
This second experiment also confirms the superiority of \cref{algorithm:split_Bregman}.
Using the same algorithmic parameters as in the previous experiment, the combined residual, \ie, the square root of \eqref{eq:residuals_alg}, is $2.1 \cdot 10^{-3}$ after $200$~iterations.
Its counterparts with again $3$ gradient steps per iteration both only achieved a combined residual in the magnitude of $5\cdot 10^{-1}$ in roughly the same time.
The numerical denoising results are presented in \cref{figure:Mesh_denoising_using_TV_Bunny}.

\section{Mesh Inpainting Problem}
\label{section:Mesh_Inpainting_Problem}

This section is devoted to mesh inpainting problems.
These problems differ from \eqref{eq:mesh_denoising_problem} in that there is no tracking data available at all.
Instead, the exact positions of a number of vertices are given, while the positions of the remaining vertices are unknown and there is no reference value known for them.
This situation often occurs during 3D scanning, \eg, at the supports or at concave regions.

When denoising shapes using \cref{algorithm:split_Bregman}, the tracking term $\cF$ adds an identity matrix to the Hessian, which counteracts the extreme nonlinearity of the rest of the augmented Lagrangian \eqref{eq:augmented_Lagrangian} \wrt to the shape.
To make use of a Newton scheme for mesh inpainting as well, we reintroduce a tracking term to the objective by a prox-like strategy, which penalizes the discrepancy between the current to the nodal positions in the previous iteration~$k$, \ie, we use
\begin{equation}
	\label{eq:mesh_inpainting_problem}
	\text{Minimize}
	\quad
	\frac{1}{2} \sum_{\vertex \in \cV_0} \abs{x_\vertex - \sequence{x}{k}_\vertex}_2^2
	+
	\tau \sum_{T\in \cT} \frac{1}{\abs{T}}
	+
	\beta \sum_{E \in \cE} \abs{E}_2 \, \abs[big]{(\logarithm{\bn_+}{\bn_-}) \cdot \bmu_+}
	.
\end{equation}
Here $\cV_0$ denotes the subset of vertices pertaining to the region to be filled in.
The positions of the remaining vertices remain fixed.

For the analogue of \eqref{eq:mesh_inpainting_problem} in \emph{image} inpainting, one often starts with black pixels as an initial guess and the connectivity between pixels is clear.
By contrast, there is no canonical initial guess for mesh inpainting.
We proceed as follows in our numerical experiments.
We mark selected subregions on a given mesh as void.
We then let \gmsh (version~3.0.6) fill the voids, using a minimal surface-like method \cite{RemacleGeuzaineCompereMarchandise:2010:1}.
This also creates the connectivity for the inpainted region.

\begin{figure}[htp]
	\centering
	\includegraphics[width = 0.32\linewidth, trim = {18cm 2cm 20cm 2cm}, clip]{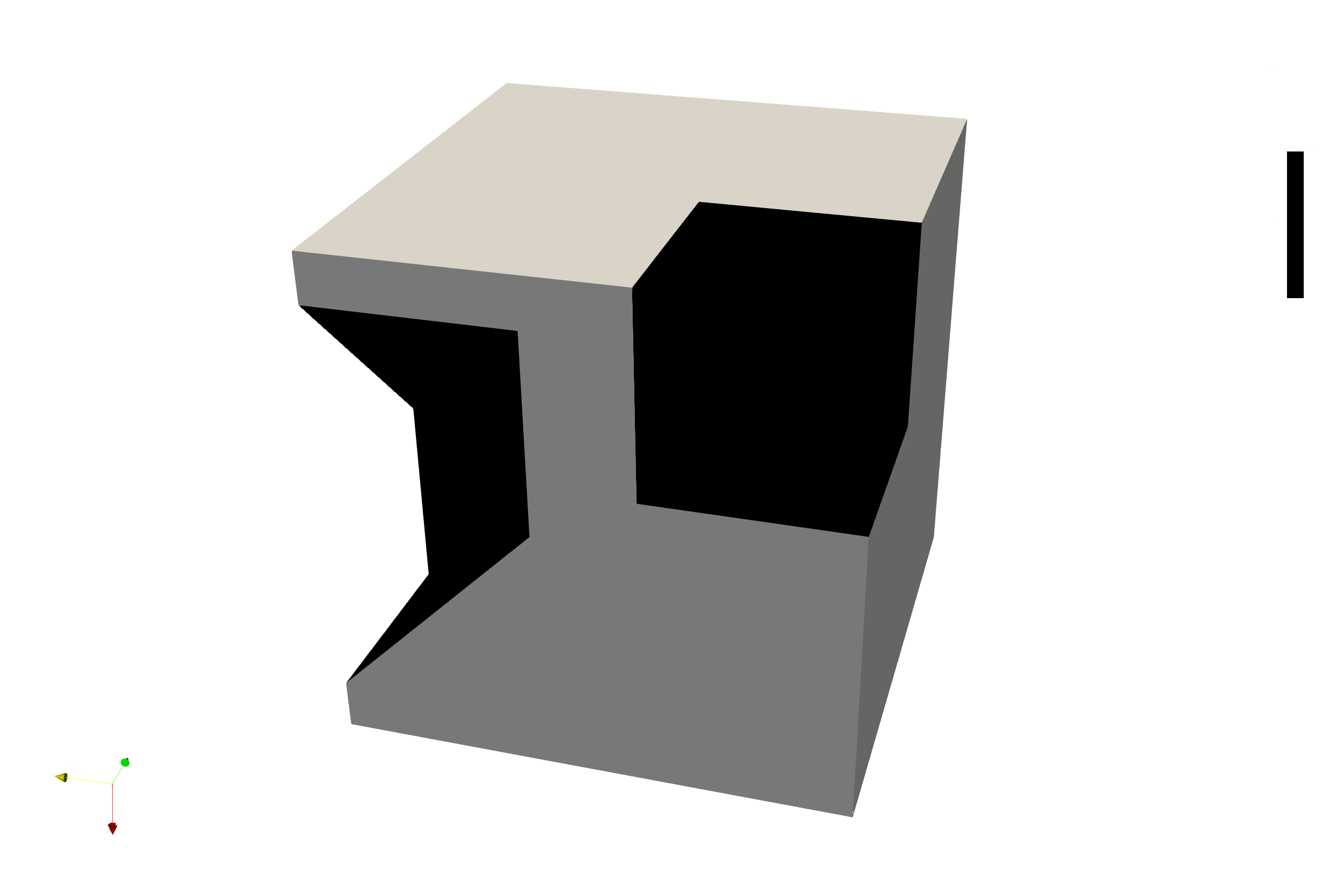}
	\hfill
	\includegraphics[width = 0.32\linewidth, trim = {18cm 2cm 20cm 2cm}, clip]{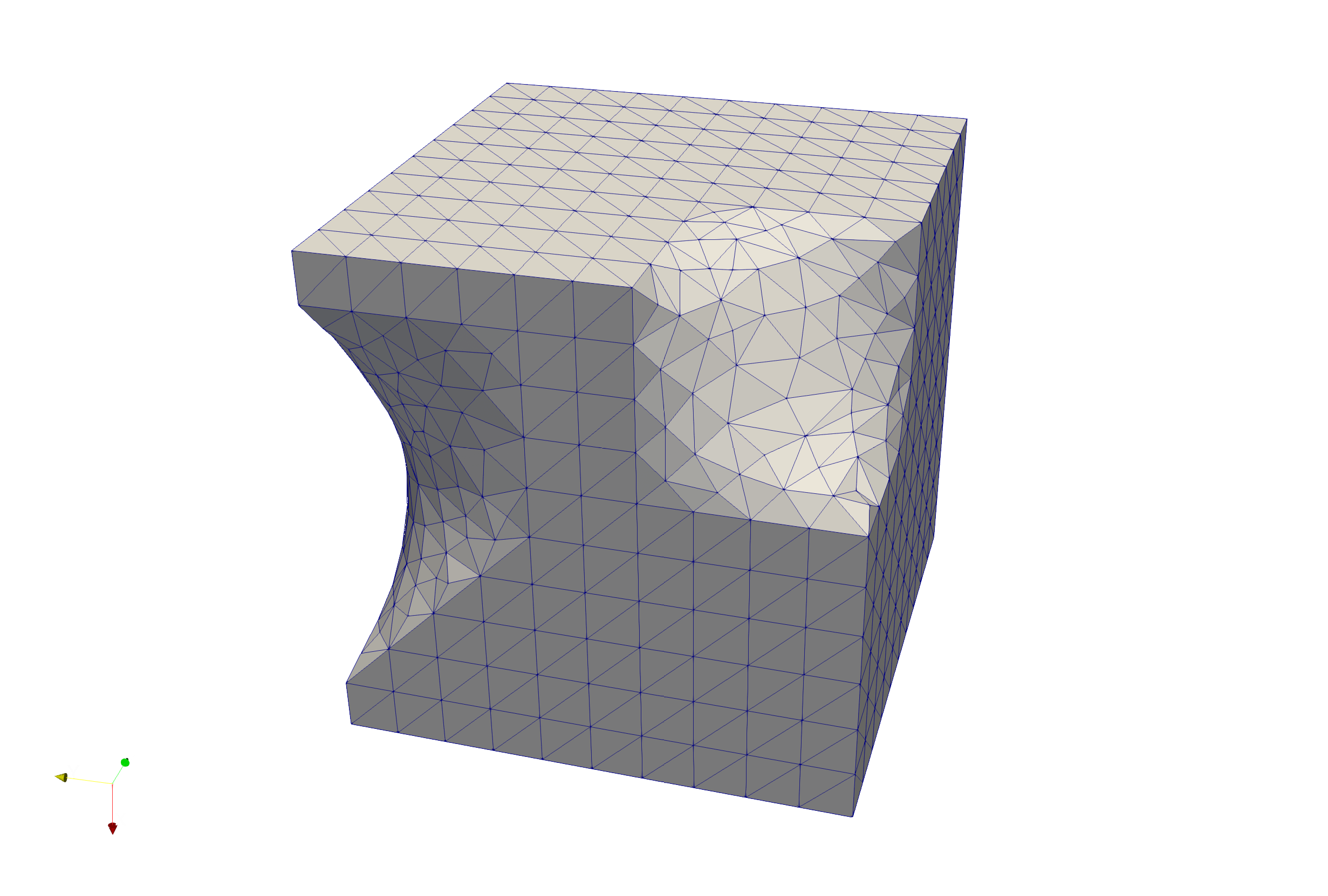}
	\hfill
	\includegraphics[width = 0.32\linewidth, trim = {18cm 2cm 20cm 2cm}, clip]{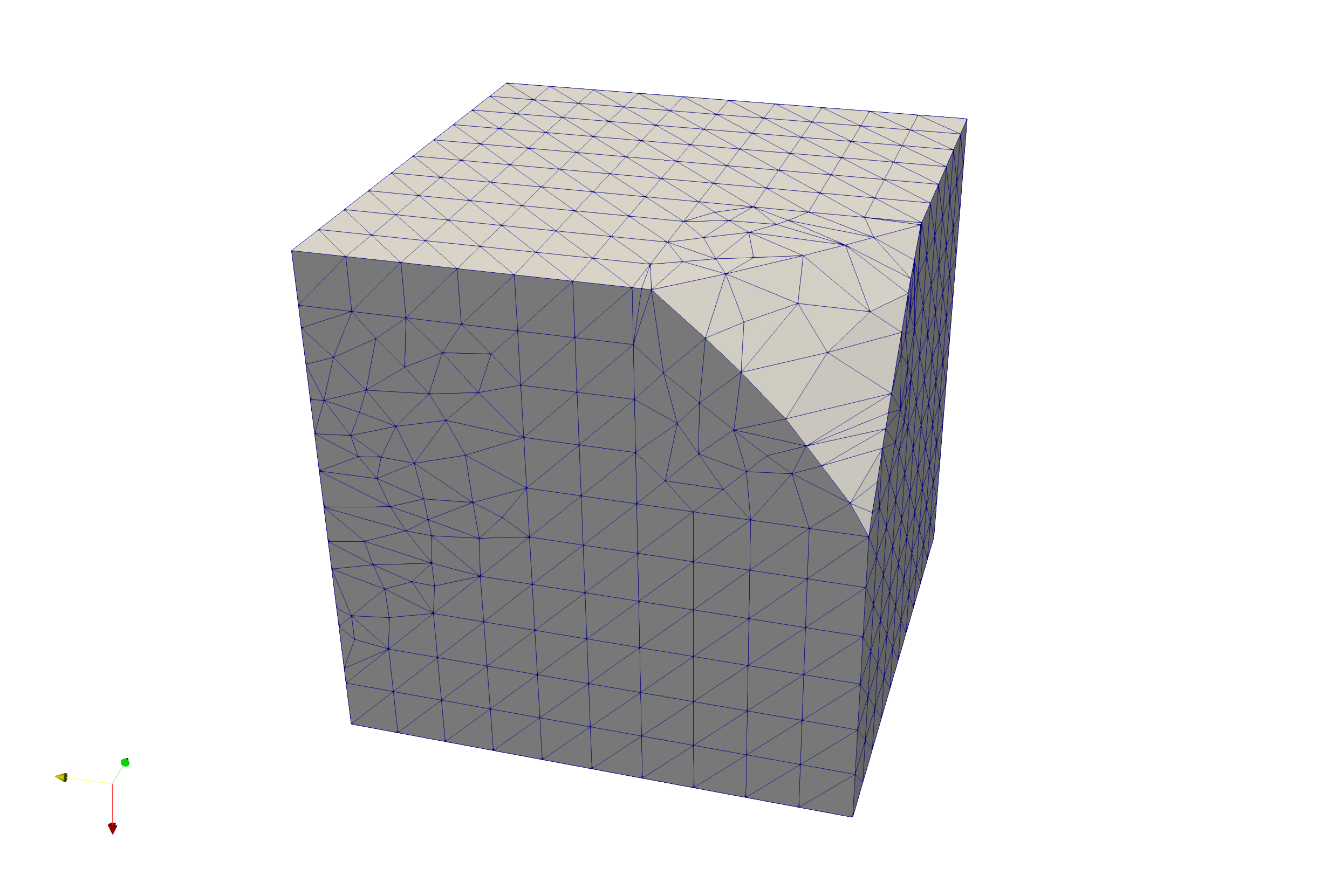}
	\caption{Mesh inpainting result using \cref{algorithm:split_Bregman} on problem \eqref{eq:mesh_denoising_problem} with $\beta = 10^{-2}$, $\tau = 10^{-8}$ and initial penalty parameter $\rho = 10^{-2}$. Original geometry (left), noisy geometry (middle) and reconstruction (right).}
	\label{figure:Cubeinpainting}
\end{figure}

We consider two test cases.
The first case is a simple unit cube mesh whose facets consist of $10 \times 10$ squares subdivided into two triangles each.
Two subregions are marked as void and an initial guess is created using \gmsh.
The inpainting results, obtained by using \cref{algorithm:split_Bregman}, are shown in \cref{figure:Cubeinpainting}.

Knowing that the ground truth is a cube, we can see that our algorithm nicely reconstructs this geometry except for the front-facing corner.
Apparently, \enquote{chopping} the corner yields a smaller value of the total variation of the normal vector than the original geometry.
This is an artifact of the model, not of the algorithm.

\begin{figure}[htb]
	\centering
	\includegraphics[width = 0.33\linewidth, trim = {10cm 0cm 3cm 0cm}, clip]{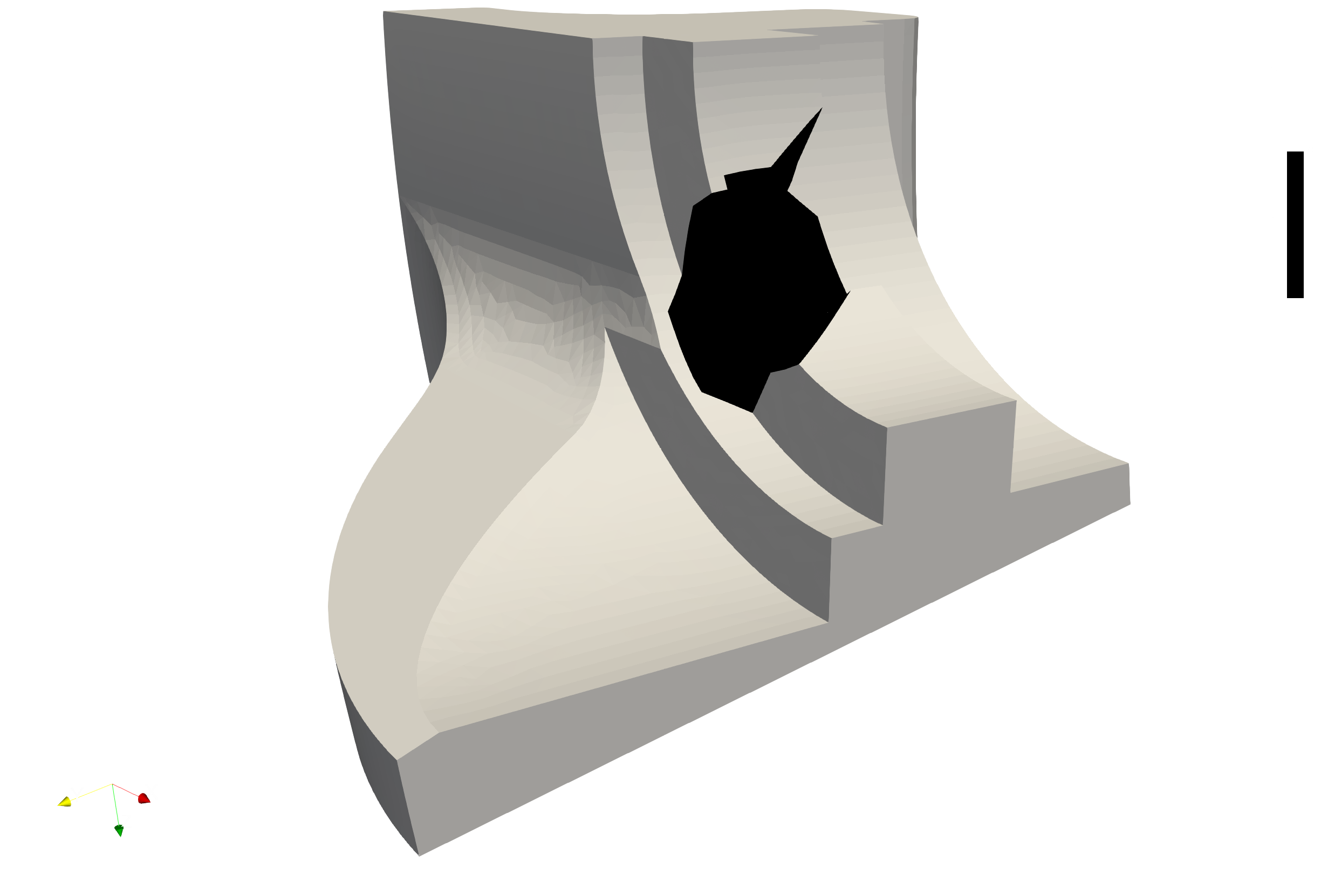}%
	\includegraphics[width = 0.33\linewidth, trim = {10cm 0cm 5cm 0cm}, clip]{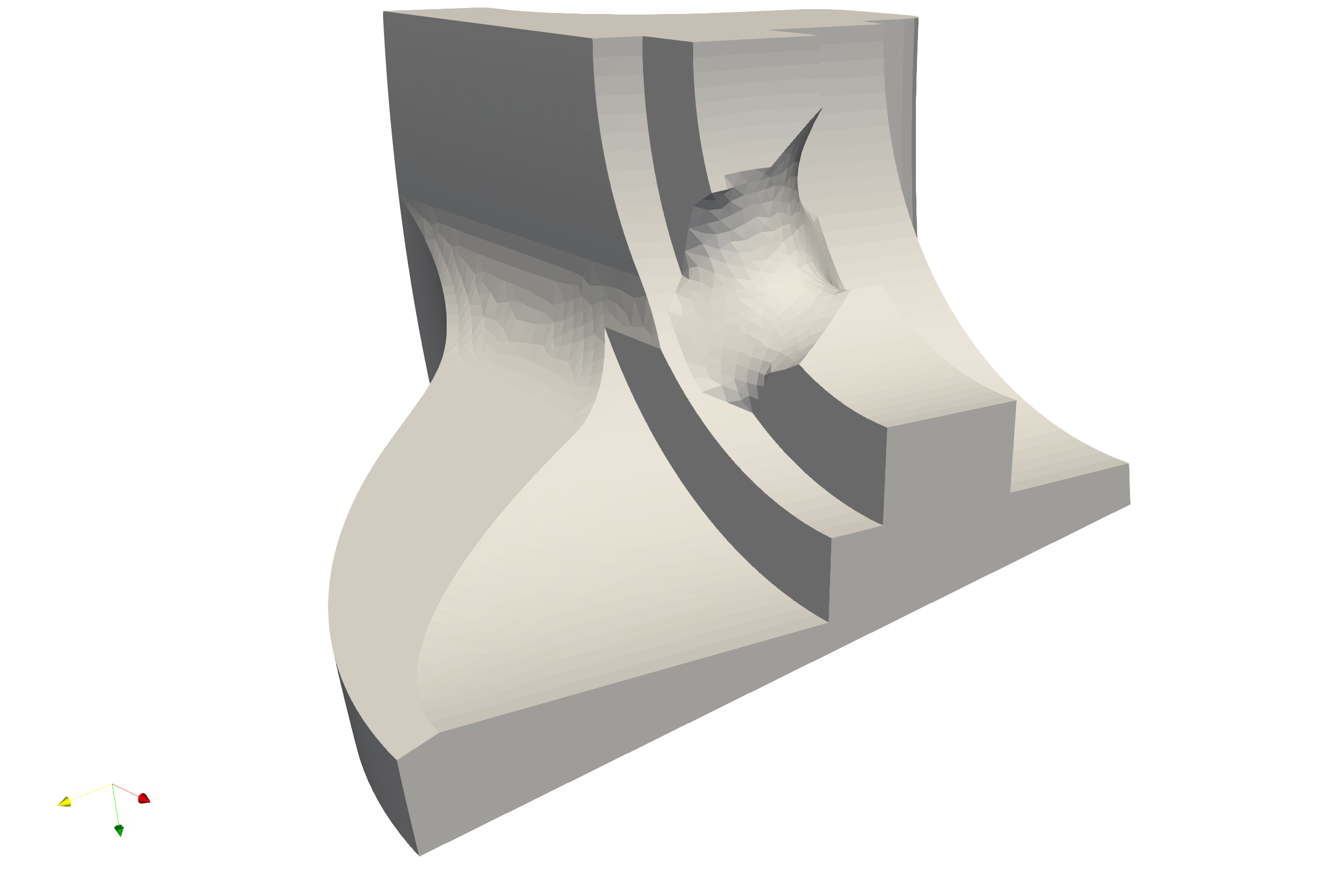}%
	\includegraphics[width = 0.33\linewidth, trim = {10cm 0cm 5cm 0cm}, clip]{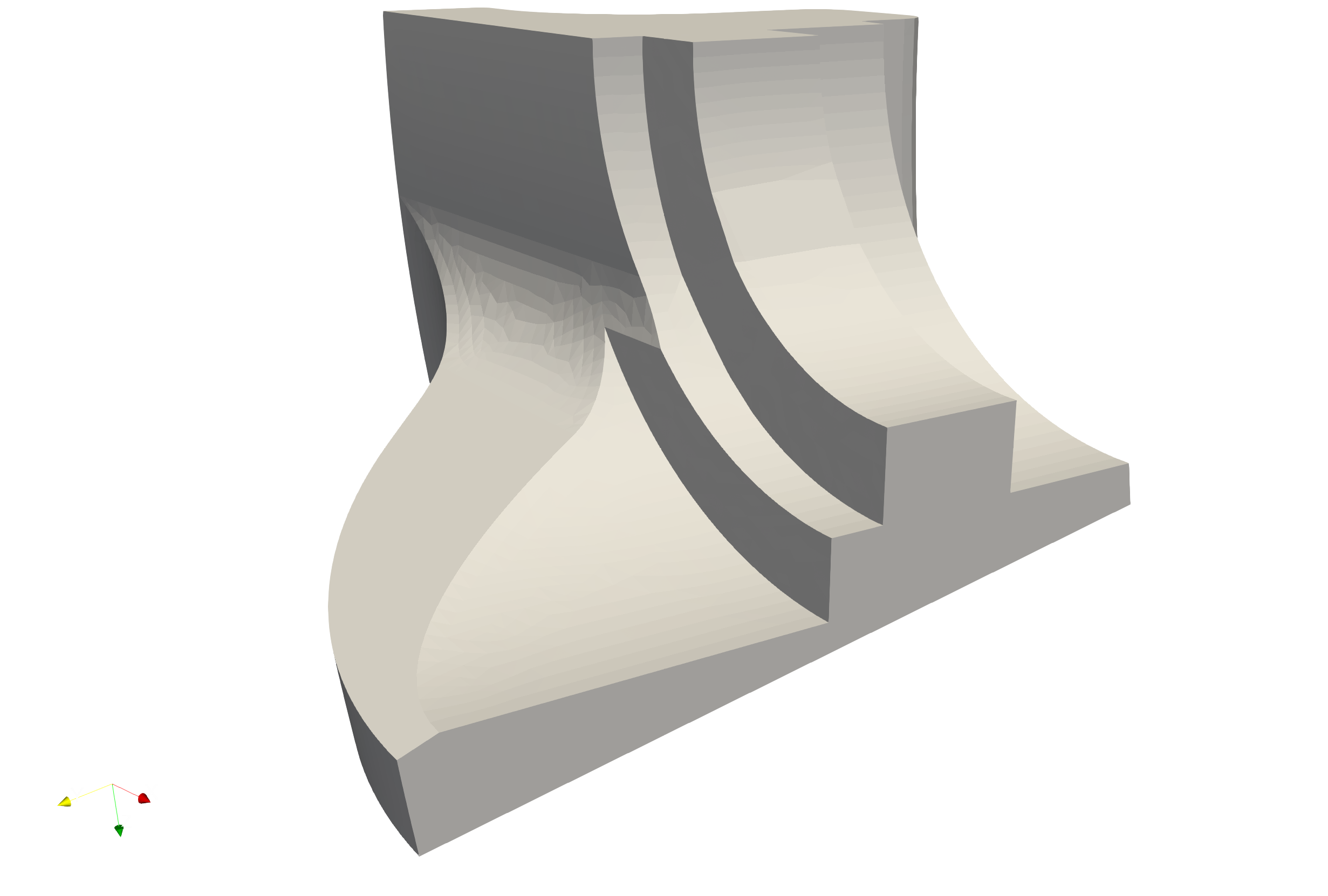}
	\caption{Mesh inpainting result using \cref{algorithm:split_Bregman} on problem \eqref{eq:mesh_denoising_problem} with $\beta = 2\cdot 10^{-2}$, $\tau = 10^{-8}$ and initial penalty parameter $\rho = 10^{-2}$. Original geometry (left), noisy geometry (middle) and reconstruction (right).}
	\label{figure:Mesh_inpatinting_using_TV}
\end{figure}

The second test case uses the more complex \emph{fandisk} geometry again.
We remove information inside a void region of significant size and with an irregular boundary.
Again, we produce an initial fill-in using \gmsh and restore the surface using \cref{algorithm:split_Bregman}.
The corresponding results are shown in \cref{figure:Mesh_inpatinting_using_TV}.

